%% file: main.tex
\documentclass[letterpaper, 10 pt, journal, twoside]{IEEEtran}

\IEEEoverridecommandlockouts                            

\usepackage[pdftex]{graphicx}
\graphicspath{{figure/}}
\usepackage{amsmath}
\usepackage{amssymb}
\usepackage{url}

\usepackage{amsthm} 
\usepackage{cite}
\usepackage{color}
\usepackage{comment}
\usepackage{fancyhdr} 

\usepackage{algorithm}
\usepackage{algorithmic}
\usepackage{bm}
\usepackage{bbm}
\usepackage{enumitem}
\usepackage{mathtools}

\newtheorem{theorem}{\textbf{Theorem}}
\newtheorem{lemma}{\textbf{Lemma}}
\newtheorem{assumption}{\textbf{Assumption}}
\newtheorem{remark}{\textbf{Remark}}

\newcommand{\bx}{\bm{x}}

\newcommand{\by}{\bm{y}}
\newcommand{\bz}{\bm{z}}
\newcommand{\bq}{\bm{q}}

\newcommand{\bs}{\bm{s}}
\newcommand{\bg}{\bm{g}}

\newcommand{\bv}{\bm{v}}

\newcommand{\bd}{\bm{d}}

\newcommand{\bh}{\bm{h}}
\newcommand{\bla}{\bm{\lambda}}
\newcommand{\bom}{\bm{\omega}}
\newcommand{\R}{\mathbb{R}}

\newcommand{\sS}{\mathcal{S}}

\newcommand{\sX}{\mathcal{X}}

\newcommand{\sA}{\mathcal{A}}
\newcommand{\sO}{\mathcal{O}}
\newcommand{\sZ}{\mathcal{Z}}
\newcommand{\sF}{\mathcal{F}}
\newcommand{\sL}{\mathcal{L}}
\newcommand{\sC}{\mathcal{C}}
\newcommand{\sJ}{\mathcal{J}}
\newcommand{\sM}{\mathcal{M}}

\newcommand{\B}{\mathbb{B}}

\newcommand{\bxi}{\bm{\xi}}

\newcommand{\bpsi}{\bm{\psi}}
\newcommand{\bmu}{\bm{\mu}}

\newtheorem{corollary}{Corollary}
\newtheorem{proposition}{Proposition}

\allowdisplaybreaks 

\title{Model-Free  Feedback Constrained Optimization Via Projected Primal-Dual Zeroth-Order Dynamics }

\author{Xin Chen, Jorge I. Poveda, Na Li
\thanks{X. Chen and N. Li are with the School of Engineering and Applied Sciences, Harvard University, USA; Emails: chen\_xin@g.harvard.edu, nali@seas.harvard.edu.
J. I. Poveda is with  the Department of Electrical, Computer, and Energy Engineering at the University of Colorado, Boulder, USA; Email: jorge.poveda@colorado.edu.}
\thanks{ 
The work was supported by NSF CNS 1947613, NSF CAREER: ECCS-1553407 and
NSF EAGER: ECCS-1839632.} 
}

\begin{document}

\maketitle

	
\input{Abstract}

\input{Introduction}

\input{Problem}

\input{Method}

\input{Performance}



\input{Simulation}

\input{Conclusion}

\input{Appendix1}

\input{Reference}

\end{document}

%% file: Abstract.tex
\begin{abstract}
In this paper, we propose a model-free feedback solution method to solve generic constrained optimization problems, without knowing the  specific formulations of the objective and constraint functions. This solution method is 
 termed projected primal-dual zeroth-order dynamics (P-PDZD) and is developed based on 
projected primal-dual gradient dynamics and extremum seeking control.
In particular, the P-PDZD method can be interpreted as a model-free   controller that autonomously drives an unknown system to the
solution of the  optimization problem using only output feedback. The P-PDZD can properly handle both the hard and asymptotic constraints, and we develop the decentralized version of P-PDZD when applied to multi-agent systems. Moreover,
 we prove that the P-PDZD achieves   semi-global practical asymptotic stability and structural robustness. 
 We then apply the decentralized P-PDZD to 
  the optimal voltage control problem in power distribution systems  with square probing signals, and the simulation results verified the optimality, robustness, and adaptivity of the P-PDZD method.

\end{abstract}


%% file: Introduction.tex
\section{Introduction} \label{sec:introduction}

This paper studies the  real-time feedback control design  that aims at autonomously steering a physical system to the solution of a constrained optimization problem. 
It   relates to the   emerging paradigm called \emph{feedback optimization} \cite{he2022model,hauswirth2021optimization} that interconnects optimization iterations in closed-loop with the optimal control of physical plants.
This type of control design has recently attracted considerable attention and can be applied in a wide range of fields \cite{hauswirth2021projected}, such as 
   electric power grids \cite{chen2021safe,dall2016optimal,tang2017real},  communication networks \cite{chen2011convergence,low1999optimization}, transportation systems \cite{como2021distributed}, etc.  However,  designing 
    such controllers for real physical systems is particularly challenging  because of two major obstacles. One is the \emph{lack of accurate system models}, as many real systems, such as power grids, are too complex to model and are subject to  unknown time-varying disturbances. The other obstacle is the need  to properly \emph{handle various constraints}, including physical laws, control saturation, artificial performance requirements, etc., which significantly complicate the problem. This paper aims to  develop model-free    feedback  solution methods for solving general constrained optimization, so that the solution methods can be interpreted as the desired feedback
    controllers that overcome these two obstacles.

To address the lack of system models or tackle intrinsically hard-to-model problems, \emph{model-free} control and optimization schemes have been widely studied. 
Instead of pre-establishing a static model from first principles or historical data, model-free approaches probe the unknown system and online learn its characteristics from real-time measurements or other feedback. Reinforcement learning (RL) \cite{sutton2018reinforcement,li2017deep} is a prominent type of model-free technique  that is
  concerned with how agents take
sequential actions in an uncertain interactive environment and
learn from the feedback to maximize the cumulative reward. Despite its huge success in games,   applying RL to
the control of physical plants is still under development
and entails many  limitations \cite{chen2022reinforcement}, such as safety and
 scalability issues,  limited  theoretical guarantee, etc.  Moreover, RL  focuses on the cumulative performance over Markov decision processes \cite{bertsekas2012dynamic}, which is beyond our scope of feedback control design that steers a  system to an optimal steady state.

 An alternative type of model-free technique is 
 \emph{extremum seeking (ES) control}  \cite{ariyur2003real,PovedaTAC17B},
 which is a classic adaptive control method and has  regained research momentum due to its theoretical advances \cite{tan2006non,tan2010extremum}. The basic idea of ES  is to perturb the system with a probing signal and estimate the gradient of the objective function solely  based   on  the output feedback, and then it drives the system to an extremum state
   via a  gradient descent flow-like structure. In this way, ES   exhibits a close connection with  zeroth-order optimization (ZO) approaches \cite{nesterov2017random}, which optimize using only function evaluations, and
     ES can be  regarded as the continuous-time version of  single-point ZO \cite{chen2021improve}. Due to its model-free feedback nature,
   ES has been applied to the control of various black-box or   hard-to-model   systems, such as   voltage control in power systems \cite{chen2021safe}, combustion control of thermal engines \cite{killingsworth2009hcci},  photovoltaic  maximum power point tracking \cite{6362193},  etc.  Hence,
     ES is particularly suitable for the model-free feedback control problem that this paper studies.

Despite the theoretical advances and practical applications,    
  one of the major limitation of existing ES schemes is that the constraints are not well enforced.  For real-world physical systems, their constraints are of different natures and can be categorized into two types: \emph{hard constraints} and \emph{asymptotic constraints} \cite{hauswirth2021optimization}. Hard constraints refer to the physical control saturation or actuator capacity limits that should be satisfied all the time, such as the generation capacity of a power plant. Asymptotic constraints refer to soft physical limits or artificial performance requirements that can be violated temporarily during transient processes but should be met
in the long-term steady state, e.g., the thermal limits of power lines and voltage limits imposed by industrial standards \cite{voltstand2020},  the comfortable temperature ranges required in building climate control, etc. 
Properly handling these two types of constraints  is essential to ensure the stability and optimality of the closed-loop system.
However, most existing ES methods are  confined to unconstrained problems, or use  barrier or penalty functions to incorporate constraints into the objective \cite{dehaan2005extremum,guay2015constrained,tan2013extremum,hazeleger2022sampled}, which may not guarantee precise constraints enforcement.
Besides, 
saddle point dynamics \cite{ye2016distributed,durr2013saddle,wang2019distributed}, projection  \cite{mills2014constrained,guay2018distributed}, and
 submanifolds \cite{durr2014extremum} have been adopted to account for the constraints in ES. But these schemes generally do not distinguish between  hard  and asymptotic constraints, or only consider one type of constraints in their problems.

In this paper, we propose a new real-time zeroth-order feedback solution method, named projected primal-dual zeroth-order dynamics (P-PDZD), for solving generic constrained optimization problems with both hard and asymptotic constraints. In particular, the proposed P-PDZD method can be  interpreted as a model-free feedback controller that
steers a physical system to the solution of the constrained optimization problem, without the need to know the specific formulations of the objective and  constraint  functions. 
By exploiting real-time system feedback, the P-PDZD method is inherently robust and adaptive to time-varying unknown disturbances.
The P-PDZD method is   developed based on projected primal-dual gradient dynamics (P-PDGD) and ES control (continuous-time ZO). The contributions and merits of the proposed P-PDZD method are   summarized below.
\begin{itemize}
    \item [1)]   \emph{(Model-Free)}.
    The P-PDZD only needs the zeroth-order information, i.e., function evaluations,  of the objective and constraint functions, while   their  mathematical formulations and  gradient information are not required.
    \item [2)]   \emph{(Constraints and Optimality)}. By using projection, the P-PDZD  guarantees 
    that the hard constraints are always satisfied. By convergence, the asymptotic constraints are respected and  optimality is achieved  when the P-PDZD reaches the steady states.
    \item [3)] \emph{(Multi-Agent)}. The P-PDZD naturally extends to the cooperative multi-agent problems, and can be implemented in a decentralized fashion,  
     where each agent performs computations only in the space of its individual decision variable and  
can preserve its own private information.
    \item [4)]  \emph{(Performance Guarantee)}. By using averaging theory and singular perturbation theory, we prove that the P-PDZD achieves semi-global practical asymptotic stability and structure robustness. Moreover, 
    we demonstrate  the optimality, robustness, and adaptivity of the P-PDZD via
    extensive
    numerical experiments. 

\end{itemize}

This paper generalizes our previous work \cite{chen2021safe} that
is  dedicated to the optimal voltage control  problem in power   systems. 
Compared to \cite{chen2021safe},
this paper considers a  generic constrained optimization problem, adopts general probing signals rather than only sinusoidal waves,  and provides the complete theoretical proof techniques and   insights, etc.

The remainder of this paper is organized as follows: Section \ref{sec:problem} introduces the problem formulation and  application examples.   Section \ref{sec:algorithm} develops the P-PDZD method and presents its multi-agent version. Section \ref{sec:analysis} analyzes the theoretical performance of P-PDZD.  
  Numerical experiments are conducted in Section \ref{sec:simulation}, and conclusions are drawn in Section \ref{sec:conclusion}.



\vspace{2pt}
\noindent\textbf{Notations.}
We use unbolded lower-case letters for scalars, and bolded lower-case letters for column vectors. 
$\mathbb{R}_+:=[0,+\infty)$ denotes the set of non-negative  real values. $\B$ denotes a closed unit ball of appropriate dimension.  
$||\cdot||$ denotes the L2-norm of a vector. $[\bx; \by] := [\bx^\top, \by^\top]^\top$ denotes the column merge of  vectors $\bx,\by$. Define the index set $[n]:=\{1,\cdots,n\}$ for a positive integer $n$.
 Denote the distance between a point $\bx\in \R^n$ and a nonempty closed set $\sX\subseteq \R^n$ as 
$
    ||\bx||_{\sX}:= \underset{\bm{\by}\in \sX}{\inf}\, ||\by-\bx||.
$ 
Define the  \emph{(point) projection} that projects a point $\bx\in\R^n$ onto the set $\sX$ as
$
    \mathrm{Proj}_{\sX}(\bx):= \underset{\by \in \sX}{\mathrm{arg\, inf}}  \, ||\by -\bx||. 
$

%% file: Problem.tex
\section{Problem Formulation and Application Examples} \label{sec:problem}

In this section, we first introduce the problem formulation and the settings of available information. Then, we present several application examples to  justify the problem settings.


\subsection{Problem Formulation} \label{sec:problem:intro}

Consider solving the constrained optimization problem:
\begin{subequations}\label{eq:gen}
\begin{align}
      \text{Obj.} \ \, &  \min_{\bx}\, f(\bx) \label{eq:gen:obj}\\
    \text{s.t.} \ \    &\ \bx\in \sX \label{eq:gen:set}\\
      &\ g_j(\bx)\leq 0, \qquad j\in[m],\label{eq:gen:ineq}
 \end{align}
\end{subequations}
where $\bx \in\R^n$ is the decision variable and $f:\R^n\to\R$ is the objective function.
$\sX\subseteq \R^n$ denotes the feasible set of $\bx$. The
 function vector $\bg:=[g_1;g_2;\cdots;g_m]:\R^n\to\R^m$ describes the inequality constraints. 
In particular, we consider that only zero-order information is available for   functions $f$ and $\bg$, which is formally stated in the following assumption.

\begin{assumption}
{(Available Information).} 
 The mathematical formulations of functions $f$ and $g_1,\cdots,g_m$ as well as their derivatives of any order are unavailable, while one can only access the   function evaluations of  functions $f$ and $g_1,\cdots,g_m$. Besides,  the feasible set $\sX$ is known.
\end{assumption}

\begin{remark}\normalfont
The motivation and rationale of the above problem setting are explained below.
\begin{itemize}
    \item [1)] The above problem    is motivated by the  feedback control design that aims at steering an unknown   plant in real time to an optimal solution of problem \eqref{eq:gen}.  We simplify the  plant dynamics model by the  static 
      input-to-output mapping functions $f(\bx)$ and $\bg(\bx)$, as we consider fast stable plant dynamics that converge immediately given any input $\bx\in\sX$ and we aim to optimize the steady-state performance.
      It is   referred to as the timescale separation property 
      with fast plant dynamics and  slow control, which 
       is commonly assumed in ES control \cite{tan2010extremum,ariyur2003real}.

\item [2)] For many  complex systems, their system models,   captured by the static 
      maps $f(\bx)$ and $\bg(\bx)$,  may be unavailable or too costly to estimate in practice. 
  Fortunately,  
     the widespread deployment of smart meters and sensors provides
 real-time  measurements of the system outputs. These measurements can be 
 interpreted as the   function evaluations of $f$ and $\bg$ and used as system feedback  to circumvent the unknown model information. 
 
   \item [3)] In the problem \eqref{eq:gen}, we
    distinguish  \emph{hard constraints} and \emph{asymptotic constraints} by the feasible set $\sX$  \eqref{eq:gen:set} and the inequalities  \eqref{eq:gen:ineq}, respectively. Thus, the feasible set $\sX$  \eqref{eq:gen:set}  should always be satisfied, while the inequalities \eqref{eq:gen:ineq}  can be violated temporarily during  transient processes but
need to be satisfied in  the steady states.
\end{itemize}

\end{remark}

\subsection{Application Examples} \label{sec:app}

Below we present three application examples to illustrate the problem setting, including optimal voltage control, building thermal control, and TCP flow control.

\vspace{3pt}
\subsubsection{Optimal Voltage Control}

Consider an electric distribution grid with the monitored node set $\sM$ and the controllable device set $\sC$. Each node $j\in\sM$  has real-time  measurement on the nodal voltage magnitude $v_j$, and 
the power injection $\bx_i$ of 
each 
 device
$i\in\sC$ can be controlled 
to maintain the voltage profiles $(v_j)_{j\in\sM}$ within an acceptable interval $[\underline{v}_j, \bar{v}_j]$. Then, the optimal voltage control \cite{chen2021model,qu2019optimal} is to achieve the optimal power injection $\bx\coloneqq (\bx_i)_{i\in\sC}$ that solves the problem  (\ref{eq:ovc}):
\begin{subequations}   \label{eq:ovc}
\begin{align}
    \text{Obj.}\ & \min_{\bx}\,  \sum_{i\in\sC} c_i(\bx_i)  \label{eq:ovc:obj} \\
\text{s.t.}\ &  
 \ \bx_i \in \mathcal{X}_i, && i\,\in\, \mathcal{C} \label{eq:ovc:x}\\
& \      \underline{v}_j\leq v_j(\bx) \leq \bar{v}_j, &&j\in \mathcal{M}. \label{eq:ovc:v}
\end{align}
\end{subequations}
Here, $c_i(\cdot)$ is the cost function and $\sX_i$ is the power capacity feasible set. A key challenge is that the  voltage function $v_j(\bx)$  is generally unknown and
hard to estimate, because it actually captures the entire power grid information, including network topology, power line parameters,  power flow equations,  uncontrollable power disturbances, etc.  To address this issue, the real-time voltage measurement $v_j$ can be leveraged as system feedback to circumvent the unknown models \cite{chen2021model,qu2019optimal}.

\vspace{3pt}
\subsubsection{Building Thermal Control}
Consider  a   building that is divided into $n$  thermal zones \cite{9146920}. For each zone $i\in[n]$,
there is an air conditioner (AC) unit operating at the power of $p_i$
with the rated  power capacity $\bar{p_i}$.  Let $T_i$ be the temperature of $i$-th zone and
$[\underline{T}, \bar{T}]$ be the   comfortable temperature interval. Then,  the building thermal control \cite{chen2020online,chen2015model,6039082}  is to optimize the AC power profile  to solve the optimization problem \eqref{eq:build}:
\begin{subequations} \label{eq:build}
\begin{alignat}{2}
      \text{Obj.} \ \, &   \min_{(p_i)_{i\in[n]}}\ \,  \sum_{i=1}^n  p_i \label{eq:build:obj} \\
    \text{s.t.} \ \  & \, p_i\in [\,0, \,\bar{p}_i\,], &&\qquad  i\in[n]       \\
    & \underline{T} \leq T_i\big((p_j)_{j\in\mathcal{S}_i}\big)\leq \bar{T}, \quad &&\qquad  i\in[n]. \label{eq:build:line}
\end{alignat}
\end{subequations}
Here, $\mathcal{S}_i \subseteq [n]$ denotes  the set of zones whose AC operations affect the temperature of $i$-th zone
via heat transfer.  The function $T_i(\cdot)$ describes  the input-to-output map from AC power to the temperature of $i$-th zone, which is complex and usually unknown in practice, as it is affected by many factors, such as the  heat capacity, 
 ambient temperature, human flow, solar irradiance, etc.
 Nevertheless, smart thermometers 
 can be deployed to measure the zone  temperature $T_i$   in real time.

\vspace{3pt}

\subsubsection{TCP Flow Control} Consider a communication network with the link set $\mathcal{L}$ and the data source set $\sS$. Denote  $r_s$ as the transmission rate of source $s\in\sS$ with the lower limit $\underline{r}_s$ and upper limit $\bar{r}_s$. The data from a source $s$ flows through a prescribed path consisting of links $\sL_s\subseteq \sL$ to its destination. The TCP flow control problem \cite{low1999optimization,magnusson2017convergence} aims to achieve  the optimal transmission rates that solve the   problem \eqref{eq:tcp}:
\begin{subequations} \label{eq:tcp}
\begin{alignat}{2}
      \text{Obj.} \ \, &   \max_{(r_s)_{s\in\sS}}\ \, \sum_{s\in \mathcal{S}}   U_s(r_s) \label{eq:tcp:obj} \\
    \text{s.t.} \ \  & \, r_s\in [\underline{r}_s, \bar{r}_s], &&\quad  s\in\mathcal{S}         \\
    &  g_l\big((r_s)_{s\in\sS_l}\big)  \geq 0 , \qquad &&\quad \, l \in \mathcal{L}, \label{eq:tcp:cap}
\end{alignat}
\end{subequations}
where $U_s: \R \to \R$ is the utility function for source $s\in\sS$, and  
$\sS_l\coloneqq \{s\, |l\in \sL_s\}$ denotes the set of sources that contain the link $l\in\sL$ in their paths. Equation \eqref{eq:tcp:cap} describes the transmission capacity constraint of each link $l\in\sL$, and function $g_l(\cdot)$ denotes the  map from the transmission rates $(r_s)_{s\in\sS_l}$ to the spare capacity of link $l$.
In practice, an accurate formulation of function $g_l(\cdot)$ may  be unavailable  due to exogenous disturbances, noises and fatigue, especially for wireless communication. The utility function $U_s(\cdot)$ may also be too complex to model.  Instead, one can observe the  realization of the utility  $U_s$ and   monitor the spare transmission capacity  $g_l$ of  each link in real time.





We note that the presented application problems \eqref{eq:build}-\eqref{eq:ovc} are simplified models  for illustrative purposes; see the cited references therein for more details. Other applications, such as resource allocation \cite{patriksson2008survey}, optimal network flow \cite{6676846}, etc., can  also fit into the problem setting described in Section \ref{sec:problem:intro}.

%% file: Method.tex
  \section{Model-Free Feedback Algorithm Design and Application to Multi-Agent Systems} \label{sec:algorithm}

In this section, we first solve  problem \eqref{eq:gen} with the projected primal-dual gradient dynamics (P-PDGD) method.
Then, we integrate ES control into P-PDGD and 
develop the projected primal-dual zeroth-order dynamics (P-PDZD) method, which is the model-free feedback solution algorithm for solving \eqref{eq:gen}. Lastly, we propose the decentralized application of P-PDZD for cooperative multi-agent systems.

We make the following two standard assumptions  on   problem \eqref{eq:gen} to render it a convex optimization  with strong duality. 

\begin{assumption} \label{ass:con_sm}
 The feasible set $\mathcal{X}$ is  nonempty, closed, and convex. The functions $f$ and $g_1,\cdots,g_m$ are convex and have locally Lipschitz gradients on $\mathcal{X}$.
\end{assumption}

\begin{assumption} \label{ass:finite}
The  problem (\ref{eq:gen}) has a finite optimum, and the Slater's conditions hold for \eqref{eq:gen}.
\end{assumption}

  These assumptions above are mainly for theoretical analysis, and the proposed methods are practically applicable to a wider range of problems that may not satisfy these assumptions. This will be validated by the numerical experiments in Section \ref{sec:simulation}.
  
\subsection{Projected Primal-Dual Gradient Dynamics}

To solve the problem \eqref{eq:gen}, we introduce the dual variable $\bla:=(\lambda_j)_{j\in[m]} \in\R_+^m$ for the inequality constraints \eqref{eq:gen:ineq}
and formulate the  saddle 
 point problem  \eqref{eq:saddle}:
\begin{align} \label{eq:saddle}
     \max_{\bla\in \R_+^m}\,  \min_{\bx\in\mathcal{X}}\, L(\bx,\bla):=f(\bx)+  \sum_{j=1}^m \lambda_j g_j(\bx),
\end{align}
where  $L(\bx,\bla)$ is the  Lagrangian function. By strong duality (Assumption \ref{ass:con_sm} and \ref{ass:finite}), the $\bx$-component of any saddle point of \eqref{eq:saddle} is an optimal solution to problem \eqref{eq:gen}. 
Denote $\bz:= [\bx;\bla]$  and 
define $\mathcal{Z}:= \mathcal{X}\times \R_+^{m}$ as the feasible set  of $\bz$ in (\ref{eq:saddle}). 

To ensure the satisfaction of the hard feasible set,  we adopt the P-PDGD \eqref{eq:cppdgd} to solve the saddle point problem \eqref{eq:saddle}: 
\begin{subequations} \label{eq:cppdgd}
\begin{align} 
         \dot{\bx} &= k_x\Big[\mathrm{Proj}_{\sX}\Big( \bx - \alpha_x\big( \nabla f(\bx) +\! \sum_{j=1}^m \lambda_j \nabla g_j(\bx) \big)   \Big)\! -\bx \Big] \label{eq:cppdgd:x} 
  \\
   \dot{\lambda}_j &  = k_\lambda \Big[\mathrm{Proj}_{\R_+}\Big( \lambda_j + \alpha_\lambda\, g_j(\bx)
     \Big) -\lambda_j\Big],\quad j\in[m], \label{eq:cppdgd:la}
\end{align}
\end{subequations}
where  $k_x,k_\lambda,\alpha_x, \alpha_\lambda$ are positive parameters. 

The P-PDGD \eqref{eq:cppdgd} is referred to as  a \emph{globally projected} (or \emph{continuous projected}) dynamical system
\cite{gao2003exponential,xia2000stability}.
As the
projection operator $\mathrm{Proj}_{\sZ}(\cdot)$ is globally Lipschitz with the Lipschitz constant $L=1$ \cite[Proposition 2.4.1]{clarke1990optimization},
one can show that
the vector field of  \eqref{eq:cppdgd} is \emph{locally Lipschitz} on $\sZ$ under Assumption \ref{ass:con_sm}. 
This is different from the widely-used discontinuous projected dynamics  \cite{nagurney2012projected,8571158,cherukuri2016asymptotic} that project the gradient flow 
onto the tangent cone of the feasible set, and thus they need  sophisticated analysis
tools for discontinuous dynamical systems. The use of continuous projection facilitates the performance analysis of the P-PDGD \eqref{eq:cppdgd} and the subsequent  P-PDZD method \eqref{eq:c:pdzd}. Nevertheless, the discontinuous version of P-PDGD, i.e.,  \eqref{eq:dppdgd}, can also be used to develop  model-free feedback solution algorithms; see Appendix \ref{app:discon} for details.

By \cite[Lemma 3]{gao2003exponential},  the solution $\bz(t)$ of the P-PDGD \eqref{eq:cppdgd} will stay within $\sZ$ for all time $t\geq t_0$ when the initial  condition $\bz(t_0)\in \sZ$. Hence,
the hard feasible set \eqref{eq:gen:set} can be always satisfied during the solution process. Moreover, the following proposition states 
the equivalence between the optimal solutions of the saddle point problem (\ref{eq:saddle}) and
the equilibrium points of P-PDGD (\ref{eq:cppdgd}). And Theorem \ref{thm:cppdgdcon} states the global asymptotic convergence of P-PDGD \eqref{eq:cppdgd}.

\begin{proposition} \label{prop:equi}
Under Assumption \ref{ass:con_sm} and \ref{ass:finite}, any equilibrium point of the P-PDGD (\ref{eq:cppdgd})   is an  optimal solution of the saddle point problem (\ref{eq:saddle}), vice versa.
\end{proposition}


 \begin{theorem} \label{thm:cppdgdcon}
(Global Asymptotic Stability of P-PDGD). 
Under  Assumption \ref{ass:con_sm} and \ref{ass:finite},  the P-PDGD (\ref{eq:cppdgd}) 
with initial condition $\bz(t_0)\in \sZ$ has a unique continuously differentiable solution
 $\bz(t):[t_0,+\infty)\to \sZ$, and the solution $\bz(t)$ 
 globally asymptotically converges to an optimal solution $\bz^*$ of the saddle point problem (\ref{eq:saddle}).
\end{theorem}

The proof of Proposition \ref{prop:equi} is provided in  Appendix \ref{app:prop1}.
The  proof of Theorem \ref{thm:cppdgdcon} mainly follows the asymptotic stability of globally projected  dynamical systems  \cite{gao2003exponential} \cite[Lemma 2.4]{bansode2019exponential}; see Appendix \ref{app:thm:cp} for a 
detailed proof.

\vspace{2mm}  
However, the  P-PDGD \eqref{eq:cppdgd} method is not implementable under  the information setting in Section \ref{sec:problem:intro}, because it needs 
 the gradient information of functions $f$ and $\bg$ in \eqref{eq:cppdgd:x}. This issue is addressed in the next subsection using ES control.
 


\subsection{Projected Primal-Dual Zeroth-Order Dynamics}\label{sec:PDZD}

To address the issue of unknown gradients, we integrate  ES control into the P-PDGD  \eqref{eq:cppdgd}
and develop the projected primal-dual zeroth-order dynamics (P-PDZD) to solve problem \eqref{eq:gen} using only zeroth-order feedback. 
 Specifically, we add
 a small periodic probing signal $\epsilon_a d(\omega_i t)$ to each element  $x_i$ of the decision variable  $\bx$, leading to
\begin{align} \label{eq:xhat}
   \hat{x}_i(t) = x_i(t) +\epsilon_a d(\omega_i t),\quad i\in[n],
\end{align}
where positive parameters $\epsilon_a$ and  $\omega_i$ denote the
amplitude\footnote{For notational simplicity, we use an identical amplitude $\epsilon_a$ for all decision variables here. In practice, different amplitude parameters can be used.} and frequency, respectively.
We make the following standard assumption on the probing signal $d(t)$  \cite{tan2008choice}.
 \begin{assumption}\label{ass:prob}
 The probing signal $d(t)$ is a periodic function with the period $T = 2\pi$ that satisfies 
  \begin{align}\label{eq:signal:const}
     \int_0^{T}\!\! d(t) dt = 0;\  \frac{1}{T} \int_0^{T}\!\! d^2(t) dt =\eta_d\! >0;  \
      \max_{t\in[0,T]}\! d(t) = 1.
 \end{align}
 \end{assumption}
Common choices of the  probing signals $d(t)$ include sinusoidal wave $d_{\sin}(t) \coloneqq \sin(t) $, square wave $d_{\mathrm{sq}}(t)$ defined as \eqref{eq:square}, triangular wave, etc. \cite{tan2008choice}
 \begin{align}\label{eq:square}
     d_{\mathrm{sq}}(t)\coloneqq \begin{cases}
     1, & t\in [2k\pi, (2k\!+\!1)\pi]\\
     -1, & t\in [(2k\!+\!1)\pi, 2(k\!+\!1)\pi]
     \end{cases},\  k\in\mathbb{Z}.
 \end{align}

In addition,
the frequency parameters $\bom:=(\omega_i)_{i\in[n]}\in \R^n$ are 
defined as
\begin{align} \label{eq:fredef}
  \qquad \omega_i =\frac{2\pi}{\epsilon_\omega}\kappa_i,\quad   i\in [n],
\end{align}
where $\epsilon_\omega$ is a small positive parameter and
$\kappa_i\neq \kappa_j$ for all $i\neq j$ are positive rational numbers.
Thus  each decision variable $x_i$ is assigned with a different frequency $\omega_i$ to be distinguished for all $i\in[n]$. Moreover, 
the frequency parameters should satisfy the orthogonality condition \eqref{eq:signal:max}:
\begin{align} \label{eq:signal:max}
    \int_0^{T_{ij}}\!\! d(\omega_i t)d(\omega_j t)\, dt =0, \quad \forall i\neq j,
\end{align}
where $T_{ij}$ is the  period of $ d(\omega_i t)d(\omega_j t)$. The sinusoidal signal $d_{\mathrm{sin}}(t)$ naturally satisfies the condition \eqref{eq:signal:max}. When the square signal $d_{\mathrm{sq}}(t)$ is used, 
  we also require $\kappa_i \neq (2k+1)\kappa_j$ for all $k\in \mathbb{Z}$ and $i\neq j$ to meet the condition \eqref{eq:signal:max}.

Denote $\bd(\bom t)\coloneqq(d(\omega_i t))_{i\in[n]}$ as the column vector that collects all   probing signals. 
Then based on  the P-PDGD (\ref{eq:cppdgd}), we develop the P-PDZD   \eqref{eq:c:pdzd} to solve the saddle point point \eqref{eq:saddle}, which only needs the  zeroth-order information (i.e., function evaluations) of $f,\bg$.
\begin{subequations} \label{eq:c:pdzd}
\begin{align}
        \dot{\bx} &= k_x\Big[\mathrm{Proj}_{\hat{\sX}\ }\big( \bx - \alpha_x \,\bxi \big) -\bx\Big] \label{eq:c:pdzd:x} 
  \\
   \dot{\bla} & =k_\lambda\Big[ \mathrm{Proj}_{\R_+^m\!}\big( \bla + \alpha_\lambda \bmu
     \big) -\bla \Big]\label{eq:c:pdzd:la}\\
     \begin{split}
          \dot{\bxi}& =\! \frac{1}{\epsilon_g} \Big[ \!-\bxi \!+\! \frac{1}{\epsilon_a\eta_d} \Big(f(\hat{\bx}(t)) \!+\! 
          \bla^\top \bg (\hat{\bx}(t)) \Big)\bd(\bom t)    \Big]
     \end{split}
    \label{eq:c:pdzd:psi}\\ 
 \dot{\bmu} &=\! \frac{1}{\epsilon_g}\Big[\! -\bmu + \bg (\hat{\bx}(t)) \Big].  \label{eq:c:pdzd:mu}
\end{align}
\end{subequations}
In \eqref{eq:c:pdzd:psi}, $\eta_d$ is a constant defined in \eqref{eq:signal:const}, which depends on the type of the probing signal. For example,
  $\eta_d = \frac{1}{2}$ when the sinusoidal signal $d_{\mathrm{sin}}$ is used, and $\eta_d = 1$ for the square signal $d_{\mathrm{sq}}$. $\epsilon_g$ is a small positive parameter.  
$
    \hat{\bx}(t) := (\hat{x}_i(t))_{i\in[n]}=   \bx(t) +\epsilon_a\bd (\bom t)
$ denotes the perturbed  decision values at time $t$.
In \eqref{eq:c:pdzd:x},
$\hat{\sX}\subseteq \sX$ is the largest closed convex shrunken feasible set  such that  $\bx + \epsilon_a\mathbb{B}\subseteq \sX$ for any $\bx\in\hat{\sX}$, and $\hat{\sX} \to \sX$ as $\epsilon_a \to 0$. Note that $\hat{\bx}(t)$ is the actual implemented action or input to the system  for evaluation, while $\bx(t)$ is an intermediate computational variable. Thus
we replace $\sX$ by $\hat{\sX}$ in \eqref{eq:c:pdzd:x} to ensure $\hat{\bx}(t)\in \sX$ all the time.

The major difference between the P-PDGD \eqref{eq:cppdgd} and the P-PDZD \eqref{eq:c:pdzd} is the introduction of new variables $\bxi:=(\xi_i)_{i\in[n]}\in\R^n$ and $\bmu:=(\mu_j)_{j\in[m]}\in\R^m$ together
with their fast dynamics \eqref{eq:c:pdzd:psi} \eqref{eq:c:pdzd:mu}. The rationale and benefit of this design are explained in Remark \ref{remark:muxi}.

\begin{remark}\label{remark:muxi}
 \normalfont Intuitively, the   introduced variables $\bxi$ and $\bmu$ can be regarded as  the real-time approximation of   $\nabla_{\bx} L$ and $\bg$,
  respectively.
 With a  sufficiently small $\epsilon_g$, the P-PDZD \eqref{eq:c:pdzd} exhibit a time-scale separation  between the slow dynamics \eqref{eq:c:pdzd:x} \eqref{eq:c:pdzd:la} and the fast dynamics \eqref{eq:c:pdzd:psi} \eqref{eq:c:pdzd:mu}. As a result,  $\xi_i$ and $\mu_j$ quickly converge to the second  terms in the right-hand sides of \eqref{eq:c:pdzd:psi} and \eqref{eq:c:pdzd:mu}, and the second terms are indeed the approximation of  $\nabla_{\bx} L(\bx,\bla)$ and $\bg(\bx)$, respectively. The advantages by introducing these fast dynamics are twofold:
  \begin{itemize}
      \item [1)]  It facilitates  theoretical analysis of the projected dynamics via averaging
theory, as the high-frequency time-varying terms associated with $\bd(\bom t)$  are moved out of the projection operator. Besides, the fast
dynamics are linear and thus can be easily handled by
singular perturbation theory.
\item [2)] The fast dynamics \eqref{eq:c:pdzd:psi} \eqref{eq:c:pdzd:mu} can be interpreted as low-pass
filters for the estimation values, which
 can diminish the oscillations and improve
the transient performance of the closed-loop system.
  \end{itemize}

\end{remark}

\begin{figure}
    \centering
      \includegraphics[scale=0.37]{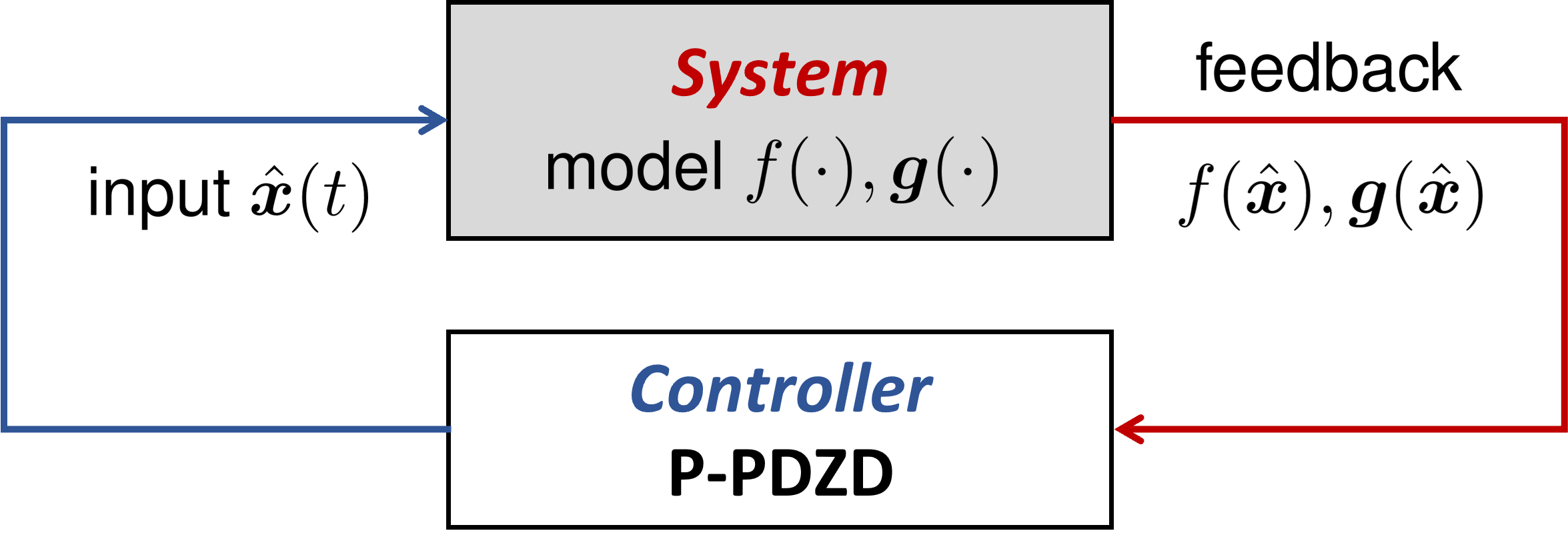}
    \caption{Implementation of the proposed P-PDZD   method \eqref{eq:c:pdzd}.}
    \label{fig:feedback}
\end{figure}

The P-PDZD \eqref{eq:c:pdzd} is the proposed optimal model-free feedback control scheme to steer the state of a black-box system
to an optimal solution of   problem \eqref{eq:gen}.
 The   implementation of the P-PDZD algorithm \eqref{eq:c:pdzd} is illustrated as 
Figure \ref{fig:feedback}.
At each time $t$, the controller feeds the input 
$\hat{\bx}(t)$ to the black-box system and receives the corresponding zeroth-order feedback  $f(\hat{\bx}(t))$ and $\bg(\hat{\bx}(t))$, which are then used to update
the input $\hat{\bx}(t)$ according to the P-PDZD \eqref{eq:c:pdzd}. 
Depending on the actual problem, the feedback can be the real-time measurements from a physical system,  the simulation outputs from a  complex simulator, or the observations from experiments.

    Although developed based on  the static optimization problem \eqref{eq:gen}, the  P-PDZD  \eqref{eq:c:pdzd}  can adapt to the dynamical system changes due to the use of real-time system feedback. 
    The adaptivity and dynamic tracking performance of P-PDZD are   demonstrated by the numerical simulations in Section \ref{sec:sim:timevary}. Moreover,  the  P-PDZD \eqref{eq:c:pdzd} directly extends to the cooperative multi-agent problems and has a decentralized version, with each agent computing and executing its own actions. This will be elaborated in the next subsection.
The theoretical analysis of the P-PDZD \eqref{eq:c:pdzd} is provided in Section \ref{sec:analysis}. Besides, the same  design idea can be applied to the discontinuous version of P-PDGD, which is presented in Appendix \ref{app:discon}.


\subsection{Decentralized P-PDZD for Multi-Agent Systems}

Consider a multi-agent system with $N$ agents.
 Each agent $i\in[N]$ is associated with a local action   $\bx_i\in \R^{n_i}$, which is subject to the feasible set
  $ \sX_i\subseteq \R^{n_i}$, and $\sum_{i=1}^N n_i =n$.
Let $\bx\coloneqq (\bx_i)_{i\in[N]}\in\R^n$ be the joint action profile of all agents.
 The goal of the agents is to cooperatively find an optimal  action profile $\bx^*$ that solves the  problem  \eqref{eq:multi}:
 \begin{subequations} \label{eq:multi}
\begin{alignat}{2}
      \text{Obj.} \ \, &   \min_{\bx}\,  \sum_{i=1}^N  f_i(\bx_i) \label{eq:multi:obj} \\
    \text{s.t.} \ \  & \, \bx_i\in \sX_i, && \qquad i\in[N]         \\
    & \, g_j\big((\bx_i)_{i\in \sS_j} \big) \leq 0, &&\qquad j\in[m],\label{eq:multi:con}
\end{alignat}
\end{subequations}
 where $\sS_j\subseteq [N]$ denotes the set of agents 
that are involved in the $j$-th constraint   described with function $g_j(\cdot)$.  Define $\sJ_i\coloneqq \{j\,| i\in \sS_j\}\subseteq [m]$ as the index set of constraints  that involve the action $\bx_i$ for each agent $i\in[N]$.
The problem \eqref{eq:multi} is a multi-agent special case of the general form \eqref{eq:gen} and is motivated by the application examples presented in Section \ref{sec:app}. Similarly, we make the following assumption on the available information for the multi-agent system.


\begin{assumption}{(Available Information for Multi-Agent System).} 
 The mathematical formulations of functions $f$ and $g_1,\cdots,g_m$ as well as their derivatives of any order are unavailable.
 Each agent $i\in [N]$ can only access the   function evaluations of   $f_i$ and $(g_j)_{j\in\sJ_i}$, 
and its own feasible set $\sX_i$.
\end{assumption}

\begin{figure}
    \centering
    \includegraphics[scale=0.38]{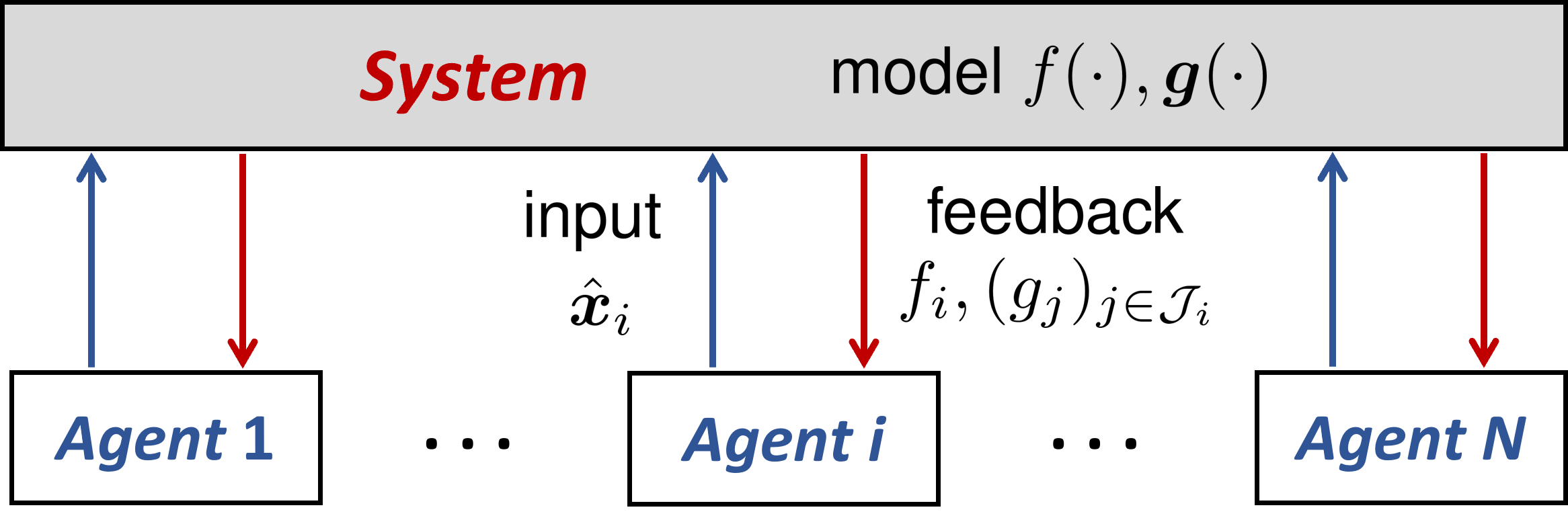}
    \caption{Implementation of decentralized  P-PDZD \eqref{eq:mzd} for  multi-agent systems.}
    \label{fig:decenPPDZD}
\end{figure}

When  the P-PDZD \eqref{eq:c:pdzd}    is applied to solve   problem \eqref{eq:multi}, it directly becomes the decentralized P-PDZD \eqref{eq:mzd}. As illustrated in Figure \ref{fig:decenPPDZD},
each agent $i\in[N]$ 
takes its own action $\hat{\bx}_i(t)\coloneqq \bx_i(t)+\epsilon_a \bd(\bom_i t)$ and receives the zeroth-order feedback $f_i(t)\coloneqq f_i(\hat{\bx}_i(t))$ and $g_j(t)\coloneqq g_j\big( (\hat{\bx}_i(t))_{i\in\sS_j}  \big)$ of all $j\in\sJ_i$. Then each agent $i$ updates its own action $\hat{\bx}_i(t)$ according to \eqref{eq:mzd}:
\begin{subequations}\label{eq:mzd}
\begin{align}
        \dot{\bx}_i &= k_x\Big[\mathrm{Proj}_{\hat{\sX}_i}\big( \bx_i - \alpha_x \bxi_i \big) -\bx_i \Big] \label{eq:mzd:x} 
  \\
      \dot{\lambda}_j & = k_\lambda\Big[ \mathrm{Proj}_{\R_+}\big( \lambda_j + \alpha_\lambda\, \mu_j
     \big) -\lambda_j\Big],   && \hspace{-15mm} j\in \sJ_i  \label{eq:mzd:la}\\
          \dot{\bxi}_i& =\! \frac{1}{\epsilon_g} \Big[ \!-\!\bxi_i +\! \frac{1}{\epsilon_a \eta_d} \big(f_i(t)\!+\!\sum_{j\in \sJ_i}\!\lambda_j g_j(t) \big)\bd(\bom_i t)    \Big]  \label{eq:mzd:psi}\\ 
 \dot{\mu}_j &= \frac{1}{\epsilon_g}\Big[ -\!\mu_j + g_j (t) \Big],  && \hspace{-15mm} j\in \sJ_i. \label{eq:mzd:mu}
\end{align}
\end{subequations}

In this way, the P-PDZD algorithm \eqref{eq:mzd} is implemented in a decentralized manner, where each agent performs computations only in the space of its own decision variable  and thus can preserve its   private information. Moreover,   the dynamics of $\lambda_j$ and $\mu_j$,  i.e., 
\eqref{eq:mzd:la} and \eqref{eq:mzd:mu}, are indeed the same for all agents $i\in\sS_j$. Hence, \eqref{eq:mzd:la} and \eqref{eq:mzd:mu} only need to be executed once and then $\lambda_j$ is broadcast to the agents $i\in \sS_j$ to update $\bx_i$ and $\bxi_i$, which can avoid repeated computations.

%% file: Performance.tex
\section{Performance  Analysis of P-PDZD} \label{sec:analysis}

In this section, we analyze the theoretical performance of the proposed P-PDZD method, including stability analysis and structural robustness to noises.

\subsection{Stability Analysis}

Denote $\bz:=[\bx;\bla]$,  $\bpsi:=[\bxi;\bmu]$,
 and $\hat{\sZ}:= \hat{\sX}\times \R^{m}_+$.
Let  $\hat{\sA}$ be the saddle point set  for the saddle point problem \eqref{eq:saddle} with  $\hat{\sX}$.
The stability properties of the P-PDZD \eqref{eq:c:pdzd} are stated as Theorem \ref{thm:c:spas}, which is proved based on averaging theory and singular perturbation theory. The detailed proof of Theorem \ref{thm:c:spas} is
provided in Appendix \ref{app:thm:escp}.

\begin{theorem}\label{thm:c:spas}
(Semi-Global Practical Asymptotic Stability). Suppose that the saddle point set $\hat{\sA}$ is compact.
Under  Assumptions \ref{ass:con_sm},  \ref{ass:finite} and \ref{ass:prob},
 there exists a class-$\mathcal{KL}$ function $\beta$ such that for any compact set $\mathcal{D}\subset \hat{Z}\times \R^{n+m}$ of initial condition $[\bz(t_0);\bpsi(t_0)]$, 
and any desired precision $\nu>0$, there exists $\epsilon_g^*>0$ such that for any $\epsilon_g\in(0,\epsilon_g^*)$, there exists $\epsilon_a^*>0$ such that for any $\epsilon_a\in(0,\epsilon_a^*)$, there exists $\epsilon_\omega^*>0$ such that for any $\epsilon_\omega\in(0,\epsilon_\omega^*)$, the solution $\bz(t)$ of the  P-PDZD (\ref{eq:c:pdzd}) satisfies 
\begin{align}\label{eq:spas}
    ||\bz(t)||_{\hat{\sA}}\leq \beta(\,||\bz(t_0)||_{\hat{\sA}},\,  t-t_0) +\nu, \ \ \forall t\geq t_0.
\end{align}
\end{theorem}

Theorem \ref{thm:c:spas} indicates that,  due to the small probing  signal $\epsilon_a\bd(\bom t)$, the solution $\bz(t)$ of  the P-PDZD (\ref{eq:c:pdzd}) will not converge to a fixed point anymore
but rather to a small $\nu$-neighborhood of $\hat{\sA}$. Nevertheless, 
by setting the parameters $(\epsilon_g, \epsilon_a, \epsilon_\omega)$ sufficiently small, one can make this precision $\nu$
 as small as desired.

The assumption of a compact saddle point set $\hat{\sA}$ in Theorem \ref{thm:c:spas} is  standard for the use of averaging theory and singular perturbation theory. For many practical applications, the feasible set
 $\hat{\sX}$  describes the physical capacity limits or control saturation and thus is naturally compact. In addition, we can 
  replace the feasible region $\R_+^m$ of the dual variable $\bla$ by the feasible box set  $[0,M_\lambda]^m$ with a sufficiently large $M_\lambda$. Thus the saddle point set $\hat{\sA} \subseteq \hat{\sX} \times [0,M_\lambda]^m$ is compact.

\subsection{Structural Robustness}

The P-PDZD method \eqref{eq:c:pdzd} heavily relies on function evaluations (or system feedback) to steer the decision  to an optimal solution of problem \eqref{eq:gen}. Then robustness is desirable to handle    small disturbances and noises that are inevitable in practice.  The  following
corollary of Theorem \ref{thm:c:spas} \cite{PovedaNaLi2019} indicates that the P-PDZD \eqref{eq:c:pdzd} is structurally robust to small additive state noise. 

\begin{corollary} \label{thm:robust}
(Structural Robustness). For any tuple of $(\epsilon_g,\epsilon_a,\epsilon_\omega)$ that induces the bound \eqref{eq:spas}, under the same conditions in Theorem \ref{thm:c:spas},
there exists $\bar{e}>0$ such that for any  additive state noise $\bm{e}(t):[t_0,+\infty)\!\to\! \R^{n}$ with $\sup_{t\geq t_0} ||\bm{e}(t)||\leq \bar{e}$, the trajectory $\bz(t)$ of the P-PDZD (\ref{eq:c:pdzd}) with the additive state  noise $\bm{e}(t)$ satisfies
\begin{align}\label{eq:robust}
    ||\bz(t)||_{\hat{\sA}}\leq \beta(||\bz(t_0)||_{\hat{\sA}},\,  t-t_0) +2\nu, \ \ \forall t\geq t_0.
\end{align}
\end{corollary}
Comparing with (\ref{eq:spas}), the P-PDZD \eqref{eq:c:pdzd} with small additive state  noise $\bm{e}$  maintains similar convergence results,
and the noise $\bm{e}$  leads to an additional precision $\nu$ term in \eqref{eq:robust}.


%% file: Simulation.tex
\section{Numerical Experiments} \label{sec:simulation}

In this section, we apply the proposed P-PDZD method to solve the optimal voltage control (OVC) problem 
 described in Section \ref{sec:app}.  We demonstrate the optimality,  robustness, and adaptivity of the P-PDZD method via
 numerical experiments.

\subsection{Experiment Setup}

\begin{figure}
    \centering
     \includegraphics[scale=0.51]{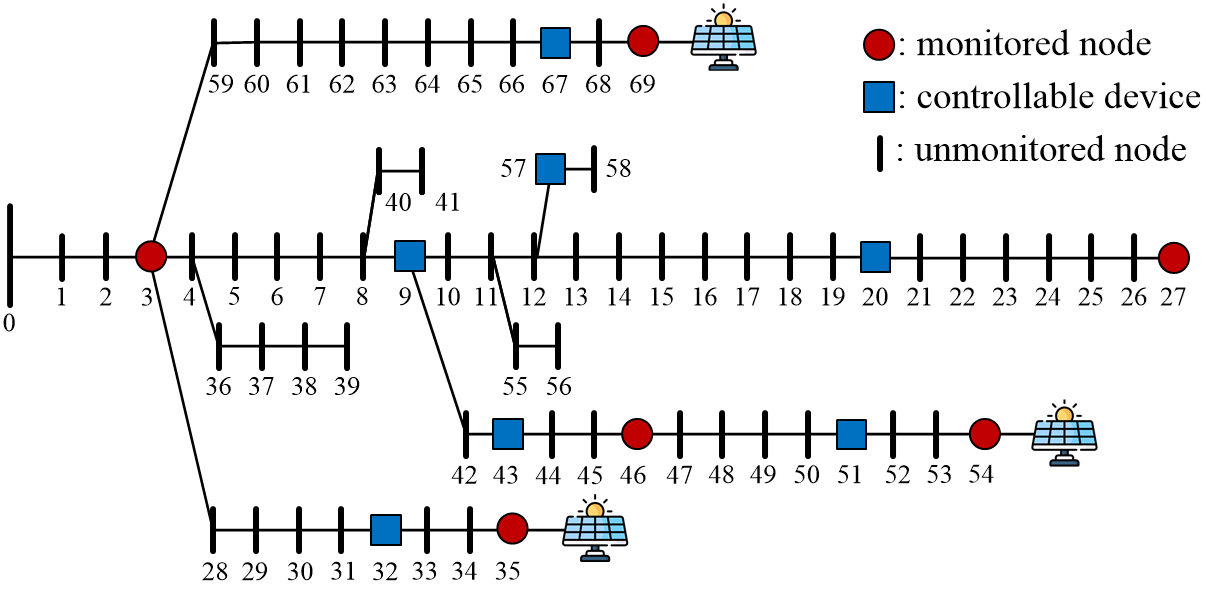}
    \caption{The modified PG\&E 69-node electric distribution grid.}
    \label{fig:PGE49}
\end{figure}

We use the  modified PG\&E 69-node electric distribution network \cite{chen2015robust}, shown as Figure \ref{fig:PGE49},  as the test system.  Three photovoltaic (PV) power plants locate at nodes 35, 54 and 69, whose  time-varying generations  are treated as unknown system disturbances that jeopardize the voltage security.
There are 7  static VAR compensators (SVCs)  located at nodes  9, 20, 32,  43, 51, 57, 67, which are the controllable devices (depicted by the set $\sC$) for voltage control under disturbances.  The reactive power output of each SVC $i\in\sC$ is 
the decision variable $x_i\in \sX_i := [-2, 2.5]$ in the unit of MVar.  We consider a known quadratic cost function $c_i(x_i) = 0.1x_i^2 $ for all $i\in\sC$. The monitored node set $\sM$ includes nodes  3, 27, 35, 46, 54 and 69, 
 which have real-time voltage measurements. 
  The voltage of node 0 (slack node) is set as 1 p.u., 
   and the lower and upper limits of voltage magnitude are set as $\underline{v}_j =0.95$ p.u. and $\bar{v}_j =1.05$ p.u. for all monitored nodes $j\in\sM$.
    The  unknown voltage function $\bv(\bx):=(v_j(\bx))_{j\in\sM}$ is 
 simulated based on the \emph{full nonconvex Distflow model} \cite{baran1989optimal2}. 
 As a result, the OVC problem \eqref{eq:ovc} does not satisfy Assumption \ref{ass:con_sm} and is nonconvex.
 The  power line impedances, nodal loads, and other parameters of the test system are provided in
  \cite{baran1989optimal}.

The OVC problem \eqref{eq:ovc}  fits into the multi-agent model \eqref{eq:multi}, and thus the decentralized P-PDZD \eqref{eq:mzd} can be applied to steer the SVC decisions $\bx:=(x_i)_{i\in\sC}$  to an optimal solution of the OVC problem \eqref{eq:ovc}. See our previous work \cite{chen2021safe} for the detailed implementation of  the decentralized P-PDZD.
Unlike the sinusoidal probing signal used in \cite{chen2021safe},   we use the \emph{square  wave} $d_{\mathrm{sq}}(t)$ \eqref{eq:square} as the probing signal, because square waves are easier to implement in practice.
For the P-PDZD algorithm,
we set $\epsilon_a =  \epsilon_\omega =  \epsilon_g = 0.025$, $\kappa_i = 1.2 + 1.5i$ for $i=1,\cdots,7$, $\alpha_x =\alpha_\lambda =0.001$, and $k_x = 50, k_\lambda =10$.

%
\subsection{Solving    Static OVC Problem} \label{sec:sim:step}
\begin{figure}
    \centering
        \includegraphics[scale=0.315]{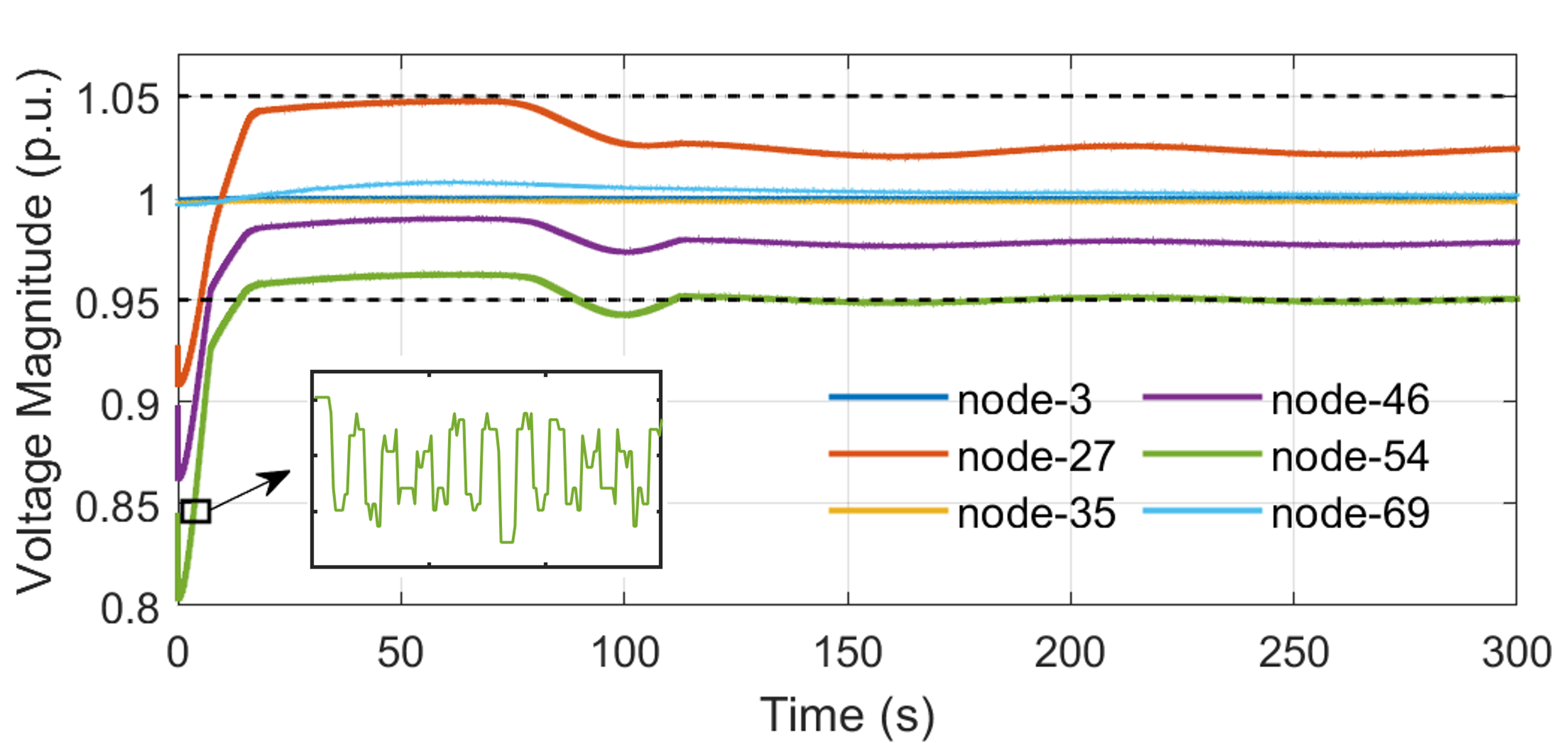}
    \caption{Voltage dynamics of the monitored nodes under  static OVC problem (black dashed lines are the upper and lower voltage limits; the curve in the small box is the zoom-in view of voltage dynamics of node-54).}
    \label{fig:stepvolt}
\end{figure}

\begin{figure}
    \centering
        \includegraphics[scale=0.315]{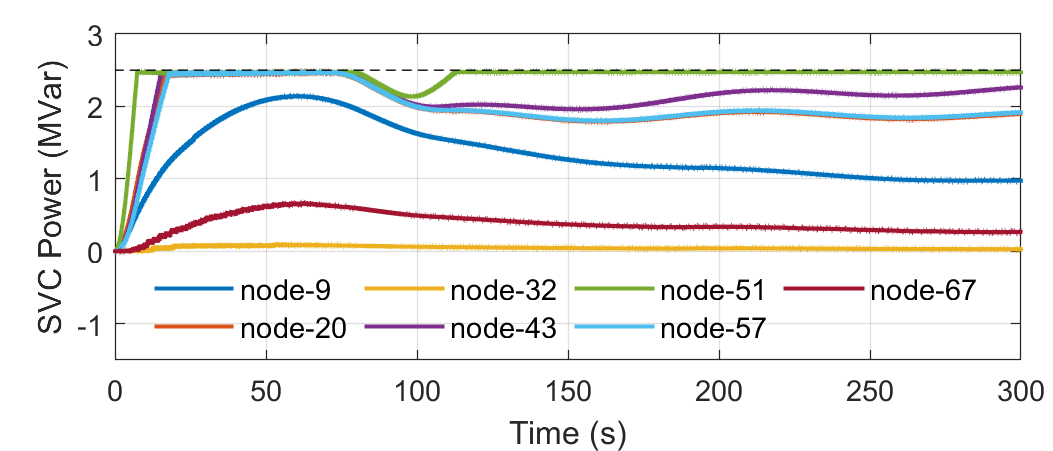}
      \caption{The  implemented reactive power outputs $\hat{\bx}(t) = \bx(t) + \epsilon_a \bd_{\mathrm{sq}}(\bom t)$ of  SVCs (the black dashed line denotes the upper power capacity limit).}
    \label{fig:power}
\end{figure}

\begin{figure}
    \centering
           \includegraphics[scale=0.315]{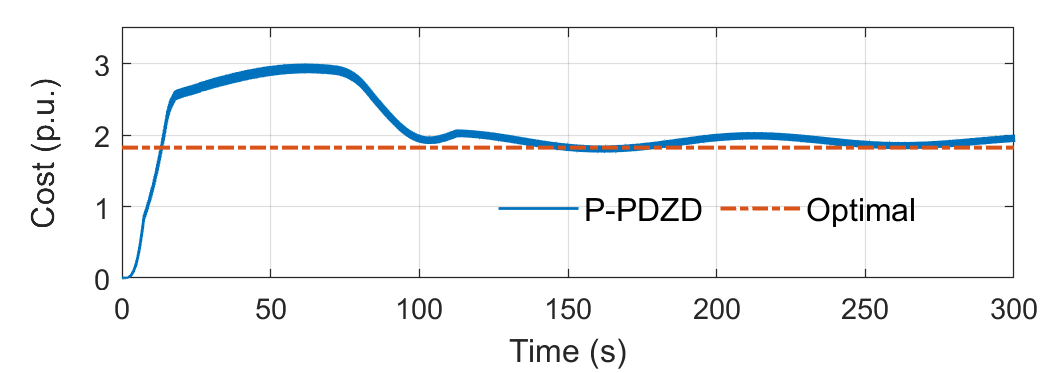}
    \caption{The dynamics of the total  control cost, i.e., the objective \eqref{eq:ovc:obj}. }
    \label{fig:cost}
\end{figure}

Consider the test scenario when the three PV plants are suddenly shut down at time $t=0$ with zero power output,
and
 voltage profiles tend to violate the lower limit due to the reduction of generation. With PV generations and loads being fixed, the OVC problem \eqref{eq:ovc} is a static optimization. We then implement the decentralized P-PDZD \eqref{eq:mzd} for real-time voltage regulation from the initial time $t=0$. The simulation results are shown as Figures \ref{fig:stepvolt}, \ref{fig:power}, and \ref{fig:cost}.
 
Figure \ref{fig:stepvolt} illustrates the voltage dynamics of all the monitored nodes. It is observed that the P-PDZD effectively brings the 
voltage profiles back to the acceptable interval $[0.95, 1.05]$ (p.u.). When zooming  in on the voltage dynamics, we can see the small-amplitude    high-frequency oscillations, which are caused by the square probing signals. Figure \ref{fig:power} shows the dynamics of SVC reactive power outputs (i.e., decision variables). It is seen that the SVC powers quickly converge  and always stay within the hard capacity limits due to the projection. The associated control cost  is shown as Figure \ref{fig:cost}, where the cost of P-PDZD converges to the optimal   value\footnote{We solve the OVC model (\ref{eq:ovc}) using the CVX package \cite{cvx} to obtain the optimal cost value (i.e., the objective \eqref{eq:ovc:obj}).}. It indicates that the P-PDZD method can steer the SVC decision to an optimal solution of the static OVC
 problem \eqref{eq:ovc}  solely
based on real-time  voltage measurement.

\subsection{Robustness to Measurement Noise}

\begin{figure}
    \centering
        \includegraphics[scale=0.315]{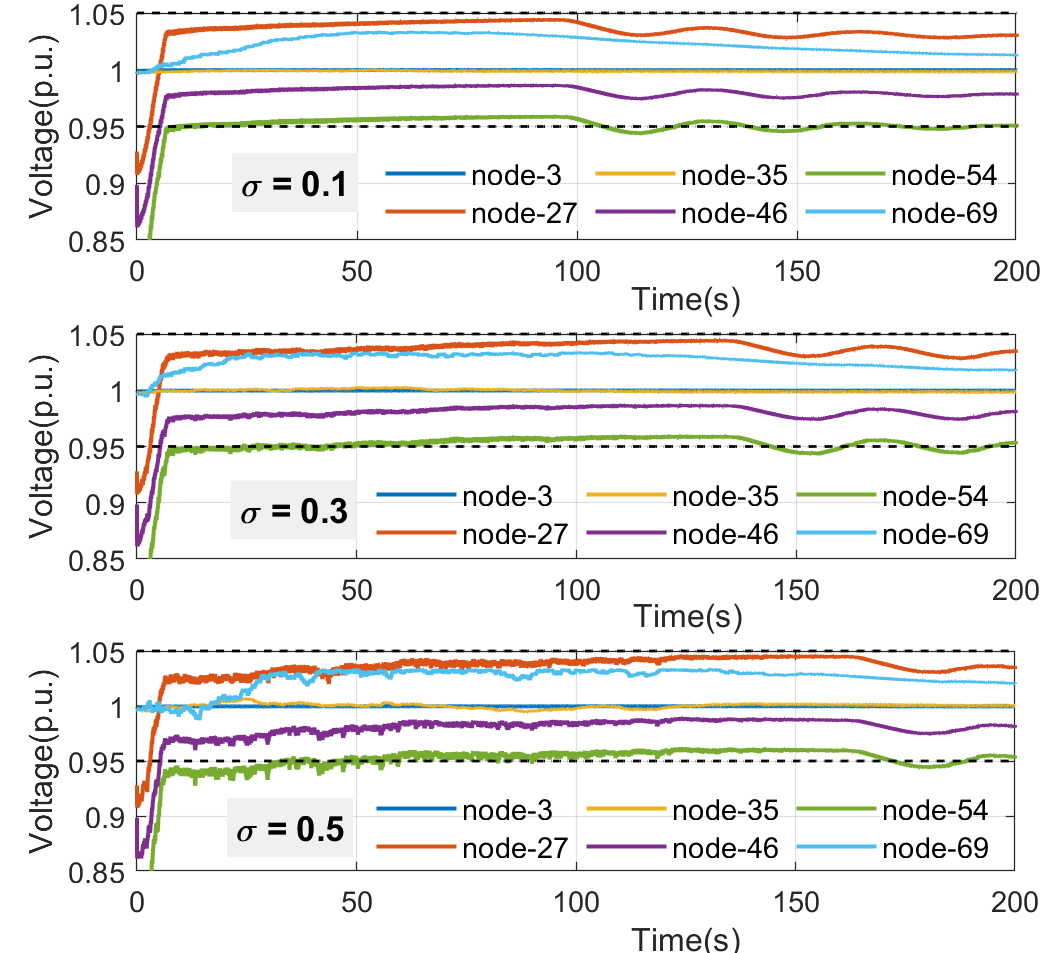}
    \caption{Voltage dynamics of the monitored nodes under different levels of  measurement noises.}
    \label{fig:noise}
\end{figure}

Measurement noises are inevitable in practice.
This subsection  considers the noisy voltage measurement $\tilde{v}^{\mathrm{mea}}_j(t)$, whose deviation from the base voltage value (1 p.u.) follows (\ref{eq:noise}):
\begin{align} \label{eq:noise}
   \tilde{v}^{\mathrm{mea}}_j(t) - 1 =  (v_j(\bx(t)) - 1)\times(1+\delta_j(t) ),
\end{align}
where $v_j(\bx(t))$ denotes the true voltage magnitude, and 
 $\delta_j$ is the perturbation ratio. Assume that $\delta_j$ is a Gaussian random variable with $\delta_j\sim \mathcal{N}(0,\sigma^2)$, which is  independent across time $t$ and  monitored nodes.  We tune the standard deviation $\sigma$ from 0 to 0.5 to simulate different levels of noises and test the performance of the P-PDZD algorithm.
 Other settings are the same as those in Section \ref{sec:sim:step}. 
   The simulation results are shown as Figure \ref{fig:noise}, while the noiseless case with $\sigma =0$ is shown in Figure \ref{fig:stepvolt}. From Figure \ref{fig:noise}, it is seen that the P-PDZD algorithm is robust to measurement noises and restores the voltage profiles to the acceptable interval in all the cases. Besides, higher levels of  noises lead to larger  oscillations in the voltage dynamics.

\subsection{Dynamic Tracking for Time-Varying OVC Problem}
\label{sec:sim:timevary}

\begin{figure}
    \centering
    \includegraphics[scale=0.315]{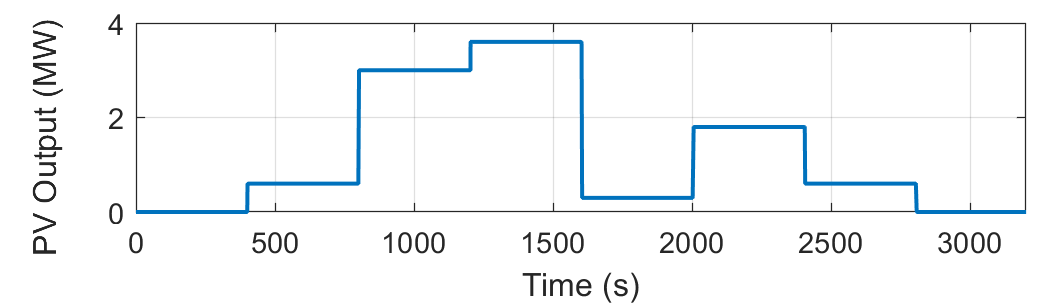}
    \caption{The time-varying total PV generation.}
    \label{fig:pvpower}
\end{figure}

\begin{figure}
    \centering
    \includegraphics[scale = 0.315]{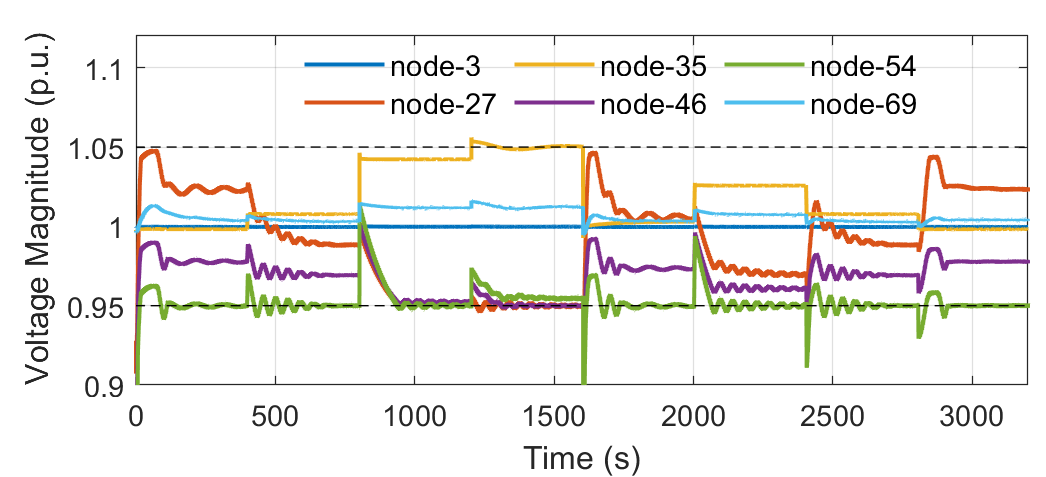}
    \caption{Voltage dynamics of the monitored buses under time-varying disturbances (black dashed lines are the upper   and lower   voltage limits).}
    \label{fig:volcompare}
\end{figure}

\begin{figure}
    \centering
    \includegraphics[scale = 0.315]{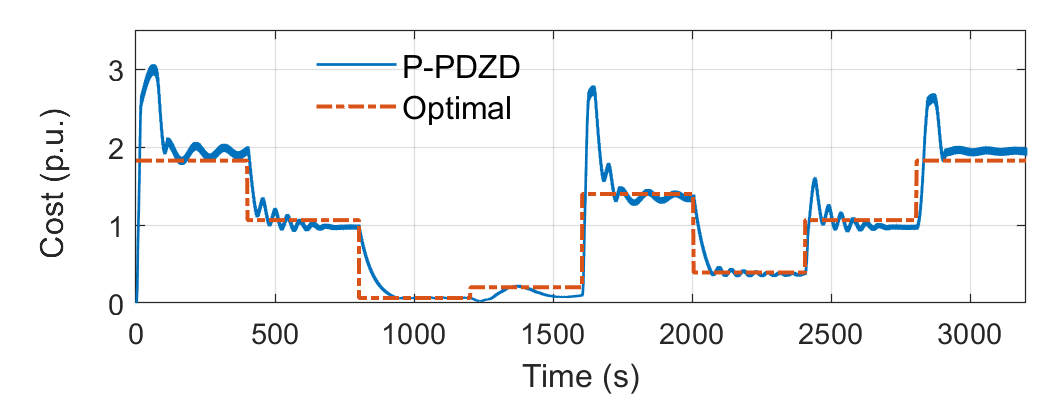}
    \caption{The control cost of P-PDZD  and the optimal cost values over time.}
    \label{fig:costconti}
\end{figure}

In practice, the power grid is a dynamic system with fluctuating loads and renewable generations. Hence, the unknown voltage function should be formulated as $\bv(\bx;\by_t)$, where $\by_t$ captures the time-varying components, 
and thus
the OVC problem \eqref{eq:ovc}
is not static but is changing over time. In this experiment, we apply the time-varying PV generation, shown as Figure \ref{fig:pvpower}, to simulate the system changes, and continuously run the P-PDZD algorithm for real-time voltage control. 
The  dynamics of the   voltage profiles and  control cost   are shown as Figures   \ref{fig:volcompare} and \ref{fig:costconti}, respectively. It is observed that   the P-PDZD algorithm  generally maintains the voltage profiles within the acceptable interval, although the voltage limits are violated very temporarily  during the transient process.  Moreover, the P-PDZD 
   keeps tracking the optimal solution of the time-varying OVC problem. This verifies that by exploiting the  real-time voltage measurement as
  system feedback, the P-PDZD can   adapt to the dynamic system  with 
   time-varying   disturbances and achieve self-optimizing performance.



%% file: Conclusion.tex
\section{Conclusion}\label{sec:conclusion}

In this paper, we propose the P-PDZD method to solve the generic constrained optimization problems 
with hard and asymptotic constraints
in a model-free feedback manner.  
Using only zeroth-order feedback or output measurements, the proposed method can be interpreted as the model-free feedback controller that autonomously drive a black-box system to the solution of the optimization problem. We prove the semi-global practical asymptotic stability and robustness of the P-PDZD and present the decentralized version of P-PDZD when applied to multi-agent problems. The numerical experiments on the optimal voltage control problem with square probing signals demonstrate the optimality, robustness, and dynamic tracking capability of the P-PDZD method. For future work, we will incorporate the plant dynamics and study the discrete-time implementation of the P-PDZD.

%% file: Appendix1.tex

\appendices

\section{Model-Free Feedback Algorithm Design Using
Discontinuous Projected Dynamics} \label{app:discon}

The discontinuous  projected primal-dual gradient dynamics (DP-PDGD) \cite{nagurney2012projected,8571158,cherukuri2016asymptotic} to solve  problem \eqref{eq:saddle} is given by
\begin{subequations} \label{eq:dppdgd}
\begin{align} 
        \dot{\bx} & =k_x\cdot\mathrm{Proj}_{T_{\sX}(\bx)}\big( -\!  \nabla f(\bx) - \sum_{j=1}^m \lambda_j \nabla g_j(\bx)  \big)  \label{eq:dppdgd:x} 
  \\
    \dot{\lambda}_j&  =k_\lambda\cdot \mathrm{Proj}_{T_{\R_+}(\lambda_j)}\big(   g_j(\bx)
     \big),\qquad j\in[m], \label{eq:dppdgd:la}
\end{align}
\end{subequations}
where $k_x,k_\lambda$ are positive parameters, and $T_{\sX}(\bx)$ denotes the  \emph{tangent cone} to  the set  $\sX$ at a point $\bx\in\sX$.

Denote $\bz:=[\bx;\bla]$ and $\sZ:= \sX\times \R_+^m$.
The DP-PDGD \eqref{eq:dppdgd}   projects the gradient flow of the Lagrangian function onto the tangent cone of the feasible set at the current point. 
When $\bz(t)$ reaches the boundary  of $\sZ$, the projection operator restricts the gradient flow such that the solution $\bz(t)$ of \eqref{eq:dppdgd} remains in $\sZ$. Hence, the DP-PDGD \eqref{eq:dppdgd} is generally  a discontinuous dynamical system \cite{cherukuri2016asymptotic}.
It is also equivalent to reformulate the DP-PDGD \eqref{eq:dppdgd} in the form of the vector projection \eqref{eq:PI-dp} \cite[Proposition 1]{brogliato2006equivalence}:
\begin{align} \label{eq:PI-dp}
    \dot{\bz} = k_z\cdot \Pi_{\sZ}(\bz, -\bh(\bz)), \ \bh(\bz):= \begin{bmatrix}
    \nabla_{\bx} L(\bx,\bla) \\
  -   \nabla_{\bla} L(\bx,\bla)
    \end{bmatrix},
\end{align}
where the  \emph{vector projection} of a direction vector $\bv $ at a point $\bz\in\sZ$ with respect to  $\sZ$ is defined as 
\begin{align}
    \Pi_{\sZ}(\bz,\bv) := \lim_{\delta \to 0^+} \frac{ \mathrm{Proj}_{\sZ}(\bz+\delta\bv )-\bz}{\delta}.
\end{align}

 Since the vector field of \eqref{eq:dppdgd} is discontinuous in general, the solution of the DP-PDGD (\ref{eq:dppdgd}) exists in the sense of an Caratheodory solution \cite{cortes2008discontinuous}. Under Assumption \ref{ass:con_sm} and \ref{ass:finite}, one can show that the DP-PDGD (\ref{eq:dppdgd}) 
with initial condition $\bz(t_0)\in \sZ$  globally asymptotically converges to an optimal solution $\bz^*$ of the saddle point problem (\ref{eq:saddle}) \cite{hauswirth2021projected, cherukuri2016asymptotic}.


%

Then we apply the same design idea proposed in Section \ref{sec:PDZD} to the DP-PDGD \eqref{eq:dppdgd} and develop the discontinuous P-PDZD (DP-PDZD) \eqref{eq:d:pdzd} for solving \eqref{eq:saddle}.
\begin{subequations} \label{eq:d:pdzd}
\begin{align} 
       \dot{\bx} & =k_x\cdot \mathrm{Proj}_{T_{\hat{\sX}}(\bx)}\big( - \!\bxi  \big)  \label{eq:d:pdzd:x} 
  \\
     \dot{\bla}& = k_\lambda\cdot\mathrm{Proj}_{T_{\R_+^m\!}(\bla)}\big(   \bmu
     \big)  \label{eq:d:pdzd:la}
     \\
    &\hspace{-2mm} \text{Equations \eqref{eq:c:pdzd:psi} \eqref{eq:c:pdzd:mu}}.
\end{align}
\end{subequations}

The DP-PDZD \eqref{eq:d:pdzd} is the same as the P-PDZD \eqref{eq:c:pdzd}, except that 
the discontinuous projection  approach is used in \eqref{eq:d:pdzd:x} and \eqref{eq:d:pdzd:la}. It indicates that
the DP-PDZD \eqref{eq:d:pdzd}  has the same properties and merits as described in the introduction section. Besides, one can also prove the semi-global practical asymptotic stability of the DP-PDZD \eqref{eq:d:pdzd} by using our proof method in Appendix \ref{app:thm:escp}, but it is more challenging because of the discontinuous dynamics, and the tools of differential inclusion and the notion of Krasovskii solutions \cite{hauswirth2021projected,cortes2008discontinuous} can be used to complete the proof.


\section{Lemmas and Proofs}

\subsection{Proof of Proposition \ref{prop:equi} }\label{app:prop1}


\noindent \textit{Proof.}
According to \cite[Theorem 3.25]{ruszczynski2011nonlinear},
the KKT conditions of the saddle point problem \eqref{eq:saddle} are given by
\begin{subequations}\label{eq:kkt}
\begin{align}
     &0\in \nabla f(\bx^*) + \sum_{j=1}^m \lambda_j^* \nabla g_j(\bx^*)+ N_{\sX}(\bx) \\
    & \bx^*\in\sX,\ g_j(\bx^*)\leq 0,\qquad\ \, \forall j\in[m]\\ 
    &\lambda_j^*\geq 0,\, \lambda_j^*g_j(\bx^*)=0,\qquad \forall j\in[m].
\end{align}
\end{subequations}

By definition,
any equilibrium point $\bz^*:=[\bx^*;\bla^*]$ of the P-PDGD \eqref{eq:cppdgd} is equivalent to satisfy \eqref{eq:ep:c}:
\begin{subequations}\label{eq:ep:c}
\begin{align}
    \bx^* &=  \mathrm{Proj}_{\sX}\Big( \bx^* \!-\! \alpha_x\big( \nabla f(\bx^*) \!+ \!\sum_{j=1}^m \lambda_j^* \nabla g_j(\bx^*) \big)   \Big) \label{eq:ep:c:x}\\
    \lambda_j^* & =  \mathrm{Proj}_{\R_+}\Big( \lambda_j^* + \alpha_\lambda\, g_j(\bx^*)
     \Big), \qquad \forall j\in[m], \label{eq:ep:c:la}
\end{align}
\end{subequations}
which is  equivalent to \eqref{eq:ep:eq}:
\begin{subequations} \label{eq:ep:eq}
\begin{align}
    &-\!\nabla f(\bx^*) \!- \!\sum_{j=1}^m \lambda_j^* \nabla g_j(\bx^*)\in N_{\sX}(\bx),\ \bx^*\in \sX\\
    &\lambda_j^*\geq 0,\ g_j(\bx^*)\leq0,\ \lambda_j^*g_j(\bx^*) = 0, \quad \forall j\in[m].
\end{align}
\end{subequations}
By comparison, equations \eqref{eq:kkt} and \eqref{eq:ep:eq} are exactly the same. Thus Proposition \ref{prop:equi} is proved. 

\subsection{Proof of Theorem \ref{thm:cppdgdcon}} \label{app:thm:cp}
We rewrite the P-PDGD \eqref{eq:cppdgd}  in a compact form\footnote{Without loss of generality, we let  $\alpha_x \!=\! \alpha_\lambda \!=\!\alpha$ and $k_x=k_\lambda=1$ for simplicity.} as
\begin{align} \label{eq:comp}
    \dot{\bz} = \mathrm{Proj}_{\sZ}(\bz - \alpha \bh(\bz)) -\bz:= \bm{f}(\bz),
\end{align}
where $\bh(\bz)$ is defined in \eqref{eq:PI-dp}. Since $\mathrm{Proj}_{\sZ}(\cdot)$ is a singleton and globally Lipschitz  with constant $L = 1$ \cite[Proposition 2.4.1]{clarke1990optimization}, the dynamics $\bm{f}(\bz)$ in \eqref{eq:comp} is locally Lipschitz  on $\sZ$ by Assumption \ref{ass:con_sm}, and thus  there exists a unique continuously differentiable 
solution $\bz(t)$ of \eqref{eq:comp} \cite[Corollary 1]{cortes2008discontinuous}. Moreover,
by \cite[Lemma 3]{gao2003exponential}, we have that $\bz(t)\in \sZ$ for all time $t\geq t_0$ whenever $\bz(t_0)\in \sZ$. 

Denote $\bz^*:=[\bx^*;\bla^*]$ as an optimal solution of \eqref{eq:saddle}.
Then, consider the following Lyapunov function  $V$:
\begin{align} \label{eq:lyap}
    V(\bz) & := \frac{1}{2}||\bz-\bz^*||^2 + L(\bx,\bla^*) -L(\bx^*,\bla) \\
    & \ =  \frac{1}{2}||\bz-\bz^*||^2 + L(\bx,\bla^*)-L(\bx^*,\bla^*)\nonumber \\
    &\quad + L(\bx^*,\bla^*)-L(\bx^*,\bla) \nonumber \geq  \frac{1}{2}||\bz-\bz^*||^2.
\end{align}
The Lie derivative   of $V$ along the P-PDGD \eqref{eq:cppdgd} is 
\begin{align}
 \mathcal{L}_{\bm{f}}V   (\bz) = \nabla_{\bz} V(\bz)^\top \bm{f}(\bz) =   (\bz-\bz^* +\bh(\bz))^\top \bm{f}(\bz).
\end{align}
One useful property is stated as  Lemma \ref{lemma:1}.
\begin{lemma}\label{lemma:1}
For any $\alpha>0$, we have
\begin{align*}
    (\bz\!-\!\bz^* +\alpha \bh(\bz))^\top \bm{f}(\bz) \leq \!-||\bm{f}(\bz)||^2 \!- \alpha (\bz\!-\!\bz^*)^\top \bh(\bz).
\end{align*}
\end{lemma}
%
The proof of Lemma \ref{lemma:1} follows \cite[Lemma 2.4]{bansode2019exponential}. For completeness,  we provide the proof  as 
the  three steps below:

\vspace{0.1cm}
\noindent 
1) We use the fact \cite{nagurney2012projected} that the projection operator satisfies 
\begin{align} \label{eq:step1}
     ( \mathrm{Proj}_{\sZ}(\bm{\gamma})-\bm{\beta})^\top (\bm{\gamma} -  \mathrm{Proj}_{\sZ}(\bm{\gamma}))\geq 0,
\end{align}
for all $ \bm{\gamma}\in\R^{|\sZ|}, \bm{\beta}\in \sZ$.

\noindent 
2) Let $\bm{\gamma} = \bz - \alpha \bh(\bz)$ and $\bm{\beta} = \bz^*$, then \eqref{eq:step1} becomes
\begin{align} \label{eq:step2}
     (\bm{f}({\bz})+\bz-\bz^*)^\top (\alpha \bh(\bz) + \bm{f}({\bz}) )\leq 0.
\end{align}

\noindent 
3) Thus we obtain Lemma \ref{lemma:1} by
\begin{align*}
  &  \, (\bz-\bz^* +\alpha \bh(\bz))^\top \bm{f}({\bz}) \\
  = &   \, (-\bm{f}({\bz}) +\bm{f}({\bz})+\bz -\bz^* +\alpha \bh(\bz))^\top \bm{f}({\bz}) \\
      = &  \, -\!||\bm{f}({\bz})||^2  + (\bm{f}({\bz})+\bz -\bz^*)^\top \bm{f}({\bz})+ \alpha \bh(\bz)^\top \bm{f}({\bz})\\
      \leq &   \, -\!||\bm{f}({\bz})||^2  - \alpha(\bm{f}({\bz})+\bz -\bz^*)^\top  \bh(\bz) + \alpha \bh(\bz)^\top \bm{f}({\bz})\\
      =&   \, -\!||\bm{f}({\bz})||^2  -  \alpha (\bz-{\bz}^*)^\top \bh(\bz),
\end{align*}
where the inequality above is because of (\ref{eq:step2}).

Using the result of Lemma \ref{lemma:1} with $\alpha=1$, we obtain 
\begin{align} \label{eq:lastineq}
    &\mathcal{L}_{\bm{f}}V   (\bz) \leq  - ||\bm{f}({\bz})||^2 - (\bz-\bz^*)^\top \bh(\bz) \nonumber\\
    = &\! -\! ||\bm{f}({\bz})||^2 \!-\!(\bx\!-\!\bx^*)^\top\nabla_{\bx}L(\bx,\!\bla)\!+\!(\bla\!-\!\bla^*)^\top\nabla_{\bla}L(\bx,\!\bla)  \nonumber \\
      \leq &\! -\! ||\bm{f}({\bz})||^2 +\! L(\bx^*,\bla) \!-\!L(\bx,\bla)  \!+\! L(\bx,\bla)\!-\!L(\bx,\bla^*) \nonumber   \\
      = &\!-\! ||\bm{f}({\bz})||^2 +\! L(\bx^*\!,\bla) \!-\!L(\bx^*\!,\bla^*) \! +\! L(\bx^*\!,\bla^*)\!-\!L(\bx,\bla^*) \nonumber \\
     \leq &\!-\! ||\bm{f}({\bz})||^2\leq 0,
\end{align}
where the second inequality follows that $L(\bx,\bla)$ is convex in $\bx$ and concave in $\bla$. 

Since $V(\bz)$ is radially unbounded and by \eqref{eq:lastineq},   the trajectory $\bz(t)$ of \eqref{eq:comp} remains bounded for all $t\geq t_0$. By LaSalle's Theorem \cite[Theorem 4.4]{khalil2002nonlinear}, we have that $\bz(t)$ converges to the largest invariant  compact subset $\mathcal{M}$ contained in $\mathcal{S}$:
\begin{align}
    \mathcal{S}:= \big\{{\bz}\in \sZ:~ \mathcal{L}_{\bm{f}}V   (\bz)=0, \, V({\bz})\leq V(\bz(t_0))      \big\}. 
\end{align}
%
When $ \mathcal{L}_{\bm{f}}V   (\bz)=0$,  we must have  $L(\bx^*,\bla) = L(\bx^*,\bla^*)$ and $L(\bx,\bla^*)=L(\bx^*,\bla^*)$ by (\ref{eq:lastineq}). Thus any point $\bz\in\mathcal{M}$ is an optimal solution of the saddle point problem (\ref{eq:saddle}).
 Lastly, following the  proof  of \cite[Theorem 15]{li2015connecting}, one can show  that $\bz(t)$ eventually converges to a fixed optimal solution $\bz^*$. 
Thus Theorem \ref{thm:cppdgdcon} is proved.

\subsection{Proof of Theorem \ref{thm:c:spas}} \label{app:thm:escp}
Denote $\bz:=[\bx;\bla]$, $\bpsi:=[\bxi;\bmu]$, and $\bs:=[\bz;\bpsi]$. 
The P-PDZD (\ref{eq:c:pdzd}) is reformulated in compact form as
\begin{align} \label{eq:original}
 \dot{\bs}=    \begin{bmatrix}
      \dot{\bz} \\   \dot{\bpsi}
    \end{bmatrix} = \begin{bmatrix}
    \bm{q}_1(\bz,\bpsi)\\
    \frac{1}{\epsilon_g}(-\bpsi + \bm{q}_2(t, \bz))
    \end{bmatrix}:= \bm{q}(t, \bs),
\end{align}
where function $\bm{q}_1(\bz,\bpsi)$ captures the dynamics \eqref{eq:c:pdzd:x} \eqref{eq:c:pdzd:la}, and function $\bq_2(t,\bz)$ is given by
\begin{align}\label{eq:g2}
    \bq_2:=\begin{bmatrix}
     \frac{1}{\epsilon_a \eta_d} \big(f(\hat{\bx}(t))+\sum_{j=1}^m \lambda_j g_j(\hat{\bx}(t)) \big)\bd(\bom t) \\
  \bg(\hat{\bx}(t)) 
    \end{bmatrix},
\end{align}
where the first part and the second part are associated with the dynamics \eqref{eq:c:pdzd:psi} of $\bxi$ and \eqref{eq:c:pdzd:mu} of $\bmu$, respectively.

We analyze the stability properties of system \eqref{eq:original} using averaging theory and singular perturbation theory,
which is divided into the following three steps.

\vspace{0.1cm}
\noindent \textbf{Step 1)} \emph{Construct a compact  set to study the behavior of system \eqref{eq:original} restricted to it.}

To apply averaging theory and singular perturbation theory, it generally requires that the considered trajectories  stay within  predefined compact sets. Without loss of generality, we consider the compact set $[(\hat{\sA}+\Delta \mathbb{B})\cap \hat{\sZ}]\times \Delta \B$ for the initial condition $\bs(t_0)$ and any desired $\Delta>0$
. Here, 
 $\hat{\sA}+\Delta \mathbb{B}$ denotes 
 the union of all sets obtained by taking a closed ball of radius $\Delta$ around each point in the saddle point set $\hat{\sA}$. 

According to Theorem \ref{thm:cppdgdcon}, there exists a class-$\mathcal{KL}$ function $\beta$ such that for any initial condition $\bz(t_0)\in \hat{\sZ}$, the trajectory $\bz(t)$ of the P-PDGD \eqref{eq:cppdgd} with the feasible set $\hat{\sX}$ satisfies 
\begin{align} \label{eq:betanew}
     ||\bz(t)||_{\hat{\sA}}\leq \beta(||\bz(t_0)||_{\hat{\sA}},\,  t-t_0) , \quad \forall t\geq t_0.
\end{align}

Without loss of generality, we assume the desired  precision $\nu\in(0,1)$. Using the $\beta$ function in \eqref{eq:betanew}, we construct the set 
\begin{align}\label{eq:Fset}
    \sF\!:=\! \Big\{\bz\!\in \!\hat{\sZ}:   ||\bz||_{\hat{\sA}} \leq \beta\big( \max_{\bv\in \hat{\sA}+\Delta \mathbb{B} } ||\bv||_{\hat{\sA}},\, 0 \big)+1  \Big\}.
\end{align}
Note that the set $\sF$ is compact under the assumption that $\hat{\sA}$ is compact.
Due to the boundedness of $\sF$, there exists a positive constant $M_1$ such that $\sF\subset M_1 \B$. Since $\bm{\ell}(\bz)$ (defined in Lemma \ref{lemma:aveg2}) is continuous by Assumption \ref{ass:con_sm}, 
 there exists a positive constant $M_2>\max\{\Delta,1\}$ such that $||\bm{\ell}(\bz)|| +1\leq M_2$ whenever $||\bz||\leq M_1$.   Denote $M_3= M_2+1$.
 We then study the behavior of system \eqref{eq:original} \textbf{restricted to evolve in the compact set $\sF\times M_3\B$}.

\vspace{0.1cm}
\noindent \textbf{Step 2)} \emph{Study the stability properties of the average system of the original system \eqref{eq:original}.}

By definition \eqref{eq:fredef}, the probing signals in system (\ref{eq:original}) are given by $d(\frac{2\pi}{\epsilon_\omega} \kappa_i t )$ for $i\in[n]$. For sufficiently small $\epsilon_\omega$, system (\ref{eq:original}), evolving in $\sF\times M_3\B$, is in  standard form  for the application of averaging theory \cite{teel2003unified}. The following Lemma \ref{lemma:aveg2} characterizes the average map for the function $\bq_2(t,\bz)$, which is proved in Appendix \ref{sec:app:lem}.

\begin{lemma}\label{lemma:aveg2}
The average of function $\bq_2(t,\bz)$ is given by
\begin{align}
    \bar{\bq}_2(\bz):=& \frac{1}{T}\!\int_0^T \!\bq_2(t,\bz)\, dt
    = \bm{\ell}(\bz) +\sO(\epsilon_a),
\end{align}
where 
$\bm{\ell}(\bz)\!:=\!\begin{bmatrix}\nabla f(\bx) + \sum_{j=1}^m \lambda_j \nabla g_j(\bx) \\
   \bg(\bx)
    \end{bmatrix}$ and $T$ is the   common  period of  $\bq_2(t, \bz)$ with fixed $\bz$.
\end{lemma}

By Lemma \ref{lemma:aveg2},
we  derive the autonomous \textbf{average system}  of 
system \eqref{eq:original} as dynamics \eqref{eq:realave} (restricted to $\sF\times M_3\B$):
\begin{align} \label{eq:realave}
    \dot{\bar{\bs}} = \!\begin{bmatrix} 
    \dot{\bar{\bz}} \\
    \dot{\bar{\bpsi}}
    \end{bmatrix} \!= \frac{1}{T}\!\int_0^T \!\!\bq(t,\bar{\bs}) dt
    \!=\!\begin{bmatrix}
    \bm{q}_1(\bar{\bz},\bar{\bpsi})\\
    \frac{1}{\epsilon_g}(-\bar{\bpsi} \!+\! \bm{\ell}(\bar{\bz}) \!+\!\sO(\epsilon_a)  )
    \end{bmatrix}
\end{align}
where $\bar{\bs}:=[\bar{\bz};\bar{\bpsi}]$ takes the same form as $\bs:=[\bz;\bpsi]$.

To analyze the average system \eqref{eq:realave}, we can first ignore the small  $\sO(\epsilon_a)$-perturbation by letting  $\epsilon_a=0$. Thus the resultant system is in the standard form for the application of singular perturbation theory \cite{Wang:12_Automatica,PovedaNaLi2019}  with the slow dynamics of $\bar{\bz}$ and fast dynamics of $\bar{\bpsi}$.
As $\epsilon_g\to 0^+$, we freeze the slow state $\bar{\bz}$, and
the \textbf{boundary layer system} of the average system \eqref{eq:realave} with $\epsilon_a=0$ in the time scale $\tau=t/\epsilon_k$ is
\begin{align}
    \frac{d \bar{\bpsi}}{d \tau} = -\,\bar{\bpsi} +\bm{\ell}(\bar{\bz}),
\end{align}
which is a linear time-invariant   system with the unique equilibrium point $\bar{\bpsi}^*=\bm{\ell}(\bar{\bz})$. As a result,
the associated \textbf{reduced system}  is derived as 
\begin{align}\label{eq:reduced}
    \dot{\bar{\bz}} = \bq_1(\bar{\bz},\, \bm{\ell}(\bar{\bz})),
\end{align}
which is exactly the P-PDGD \eqref{eq:cppdgd}. By Theorem \ref{thm:cppdgdcon} and \cite[Theorem 2]{Wang:12_Automatica}, it follows that as $\epsilon_g\to 0^+$, the set $\hat{\sA}\times {M}_3\mathbb{B}$ is semi-globally practically asymptotically stable (\textbf{SGPAS}) for the average system \eqref{eq:realave} with  $\epsilon_a=0$.
Then
by the structural robustness results
for ordinary differential equations with continuous right-hand sides \cite[Proposition A.1]{PovedaNaLi2019},  the set $\hat{\sA}\times {M}_3\mathbb{B}$ is also SGPAS for the average system \eqref{eq:realave}  as $(\epsilon_g,\epsilon_a)\to 0^+$, which is stated as Lemma \ref{lemma:y:sgpsa}.

\begin{lemma}\label{lemma:y:sgpsa} 
Given the precision $\nu$,
there exists $\epsilon_g^*>0$ such that for any $\epsilon_g\in(0,\epsilon_g^*)$, there exists $\epsilon_a^*>0$ such that for any $ \epsilon_a\in (0,\epsilon_a^*)$, 
every solution $\bar{\bs}(t)$ of the average system \eqref{eq:realave} (restricted in $\sF\times M_3\B$) with initial condition $\bar{\bs}(t_0)\in [(\hat{\sA}\!+\!\Delta \mathbb{B})\cap \hat{\sZ}]\!\times\! \Delta \B$ satisfies 
\begin{align} \label{eq:y1con}
     ||\bar{\bz}(t)||_{\hat{\sA}}\leq \beta(||\bar{\bz}(t_0)||_{\hat{\sA}},\,  t-t_0) +\frac{\nu}{4},\quad \forall t\in \mathrm{dom}(\bar{\bs}).
\end{align}
\end{lemma}
Since the average system \eqref{eq:realave} is restricted in $\sF\times M_3\B$, we have $||\bar{\bpsi}(t)||_{M_3\B} = 0$ for all $t\in \mathrm{dom}(\bar{\bs})$, which implies that $||\bar{\bs}(t)||_{\hat{\sA}\times M_3\B} =  ||\bar{\bz}(t)||_{\hat{\sA}}$ for all $t\in \mathrm{dom}(\bar{\bs})$. Hence, it follows that for all $ t\in \mathrm{dom}(\bar{\bs})$,
\begin{align*}
   ||\bar{\bs}(t)||_{\hat{\sA}\times M_3\B}\leq \beta(||\bar{\bs}(t_0)||_{\hat{\sA}\times M_3\B},\,  t-t_0) +\frac{\nu}{4}.
\end{align*}
Next we show the completeness of solutions of the average  system \eqref{eq:realave} by leveraging Lemma \ref{lemma:forward}, which follows a special case of \cite[Lemma 5]{NunoShamma_2020}.
\begin{lemma} \label{lemma:forward}
Let $M_2>0$ be given and $\bm{u}:\R_+ \to M_2\B$. Then, for any $k>0$, the set $M_2\B$ is forward invariant under the dynamics $\dot{\bpsi} = k(-\bpsi +\bm{u}(t))$.
\end{lemma}

Specifically, under the initial condition $\bar{\bs}(t_0)\in [(\hat{\sA}\!+\!\Delta \mathbb{B})\cap \hat{\sZ}]\!\times\! \Delta \B$, by \eqref{eq:y1con}, the trajectory $\bar{\bz}(t)$ of \eqref{eq:realave} satisfies $\bar{\bz}(t)\in \mathrm{int}(\sF)$ for all $t\in \mathrm{dom}(\bar{\bs})$. It implies that  $||\bz(t)||\leq M_1$ and thus $||  \bm{\ell}(\bar{\bz}(t)) \!+\!\sO(\epsilon_a)  ) ||< M_2$ for all $t\in \mathrm{dom}(\bar{\bs})$, where we take $||\sO(\epsilon_a)||<1$ for all $\epsilon_a\in (0,\epsilon_a^*)$ without loss of generality. According to Lemma \ref{lemma:forward}, $\bar{\bpsi}(t) \in M_2\B\subset \mathrm{int}(M_3\B)$ for all $t\geq t_0$. Hence,  every trajectory $\bar{\bs}(t)$ of \eqref{eq:realave} satisfies 
\begin{align}\label{eq:unbs}
    \bar{\bs}(t)\in \mathrm{int}(\sF\times M_3\B), \quad \forall t\geq t_0,
\end{align}
 and thus it has an unbounded time domain, i.e., $\mathrm{dom}(\bar{\bs})= [t_0,+\infty)$.

\vspace{0.1cm}
\noindent \textbf{Step 3)} \emph{Link the stability property of the average system \eqref{eq:realave} to the stability property of the original system \eqref{eq:original}.}

Since  the set $\hat{\sA}\times M_3\B$ is SGPAS for the average system \eqref{eq:realave} (restricted in $\sF\times M_3\B$) as $(\epsilon_g,\epsilon_a)\to 0^+$, by averaging theory for perturbed systems \cite[Theorem 7]{PovedaNaLi2019}, it directly obtains that for each pair of $(\epsilon_g,\epsilon_a)$ inducing the bound \eqref{eq:y1con}, there exists $\epsilon_\omega^*>0$ such that for any $\epsilon_\omega\in(0,\epsilon_\omega^*)$, the solution $\bs(t)$ of the original system \eqref{eq:original} (restricted to $\sF\times M_3\B$) satisfies 
\begin{align}
    ||\bs(t)||_{\hat{\sA}\times M_3\B}   \leq \beta( ||\bs(t_0)||_{\hat{\sA}\times M_3\B},\,  t-t_0) +\nu,
\end{align}
for all $t\in \mathrm{dom}(\bs)$. Since  $ ||\bz(t)||_{\hat{\sA}}= ||\bs(t)||_{\hat{\sA}\times M_3\B}$ for all $t\in \mathrm{dom}(\bs)$, we obtain the bound  \eqref{eq:spas}.  The only task left is to show the completeness of solutions of the original  system \eqref{eq:original}, i.e., proving  $\mathrm{dom}(\bs) = [t_0,+\infty)$. This can be done by   the following two lemmas.
See Appendix \ref{app:omega} for the
 proof of Lemma \ref{lemma:omega}, and Lemma \ref{lemma:close} follows \cite[Theorem 1]{Wang:12_Automatica}.

\begin{lemma}\label{lemma:omega}
There exists $\epsilon_g^*>0$ such that for any $\epsilon_g\in(0,\epsilon_g^*)$, there exists $\epsilon_a^*>0$ such that for any $ \epsilon_a\in (0,\epsilon_a^*)$,  there exists a compact Omega-limit set\footnote{See \cite[Definition 6.23]{Goebel:12} for the  notion of ``Omega-limit set of a set".} $\Omega(\sF\times M_3\B)$ that is uniformly globally asymptotically stable for the average system \eqref{eq:realave} restricted to $\sF\times M_3\B$.
\end{lemma}

\begin{lemma}\label{lemma:close}
Let $(\epsilon_g,\epsilon_a)>0$ take sufficiently small values such that  Lemmas \ref{lemma:y:sgpsa} and \ref{lemma:omega} hold. Then, for each $\Tilde{T},\delta>0$, there exists $\epsilon_\omega^*>0$ such that for all $\epsilon_\omega\in (0,\epsilon_\omega^*)$ and all solutions $\bs$ of the original system \eqref{eq:original} (restricted to $\sF\times M_3\B$), there exists a solution $\bar{\bs}$ of the average system \eqref{eq:realave}  (restricted to $\sF\times M_3\B$) such that for all $t\in \mathrm{dom}(\bs)\cap \mathrm{dom}(\bar{\bs})$, 
\begin{align}
    \sup_{t\in[0,\tilde{T}]}||\bz(t)-\bar{\bz}(t)||\leq \delta,~\sup_{t\in[0,\tilde{T}]}||\bpsi(t)-\bar{\bpsi}(t)||\leq \delta.
\end{align}
\end{lemma}

By Lemma \ref{lemma:omega}, there exists a $\mathcal{KL}$-class function $\tilde{\beta}$ such that every solution of the average system \eqref{eq:realave} satisfies 
\begin{align}
    ||\bar{\bs}(t)||_{\Omega( \sF\times M_3\B )}\leq \tilde{\beta}( ||\bar{\bs}(t_0)||_{\Omega( \sF\times M_3\B )},t-t_0),  
\end{align}
for all $t\in\mathrm{dom}(\bar{\bs})= [t_0,+\infty)$. According to 
 \cite[Theorem 2]{Wang:12_Automatica}, there exists $\epsilon_\omega^*$ such that for all $\epsilon_\omega\in(0,\epsilon_\omega^*)$, every solution of the original system \eqref{eq:original} (restricted to $\sF\times M_3\B $) satisfies 
 \begin{align}
    ||{\bs}(t)||_{\Omega( \sF\times M_3\B )}\leq \tilde{\beta}( ||{\bs}(t_0)||_{\Omega( \sF\times M_3\B )},t-t_0)+\frac{\nu}{3},  
\end{align}
for all $t\in \mathrm{dom}({\bs})$. Thus, there exists time $T_3>0$ such that 
\begin{align}\label{eq:someb}
    ||{\bs}(t)||_{\Omega( \sF\times M_3\B )}\leq \frac{\nu}{2}, \ \forall  t\in [T_3,+\infty)\cap \mathrm{dom}({\bs}).
\end{align}
Lemma \ref{lemma:close} indicates that all solutions of the original system \eqref{eq:original} remains $\delta$-close to some solution of the average system \eqref{eq:realave} on a compact time domain. Moreover, \eqref{eq:unbs} indicates that every solution of \eqref{eq:realave} stays within $ \mathrm{int}(\sF\times M_3\B)$ for all $t\geq t_0$. Thus, by applying Lemma \ref{lemma:close} with $\tilde{T} = T_3+1$, there exists $\epsilon_\omega^*>0$ such that for all $\epsilon_\omega\in (0,\epsilon_\omega^*)$ we have
\begin{align}
    \bs(t) \in  \mathrm{int}(\sF\times M_3\B), \quad \forall t\in [t_0,T_3+1].
\end{align}
In addition, by \eqref{eq:someb}, we also have $ \bs(t) \in  \mathrm{int}(\sF\times M_3\B)$ for all $t\geq T_3$. Therefore, every solution $\bs(t)$ of the original system has an unbounded time domain, i.e., $\mathrm{dom}(\bs)=[t_0,+\infty)$. Thus Theorem \ref{thm:c:spas} is proved.

\subsection{Proof of Lemma \ref{lemma:aveg2}} \label{sec:app:lem}

We first consider the  integration on the first part of $\bq_2(t,\bz)$. By the Taylor expansion of $f(\cdot)$, we have  $(\forall i\in[n])$ 
\begin{align*}
 & \frac{1}{T}\int_{0}^T \frac{1}{\epsilon_a \eta_d}  f(\bx+ \epsilon_ad(\bom t )) d(\omega_i t )\, dt \\
 = &  \frac{1}{T}\int_{0}^T\! 
\frac{1}{\epsilon_a \eta_d}  \big[f(\bx)\!+\! \epsilon_a  \nabla f(\bx)^\top d (\bom t)
\! +\!  \sO(\epsilon_a^2)\big] d(\omega_i t ) dt \\
= & \frac{1}{T}\int_{0}^T \frac{1}{\eta_d} \sum_{j=1}^n \big[\frac{\partial f(\bx)}{\partial x_j} d(\omega_j t)d(\omega_i t )\big]\, dt +\sO(\epsilon_a)\\
= & \frac{\partial f(\bx)}{\partial x_i} \frac{1}{ \eta_d T}\int_{0}^T \!   d(\omega_i t )^2 \, dt +\sO(\epsilon_a)= \frac{\partial f(\bx)}{\partial x_i}  +\sO(\epsilon_a).
\end{align*}
The third equality above is due to \eqref{eq:signal:max}.
Similarly, we have $(\forall j\in[m], i\in[n])$
\begin{align*}
  \frac{1}{T}\!\!\int_{0}^T\! \!\frac{1}{\epsilon_a \eta_d}   \lambda_jg_j(\bx\!+\! \epsilon_a\bd(\bom t )) d(\omega_i t ) dt\!= \! \lambda_j \frac{\partial g_j(\bx)}{\partial x_i} \!+\!\sO(\epsilon_a).
\end{align*}

As for the integration on the second part of $\bq_2(t,\bz)$, i.e., $\bg(\hat{\bx}(t))$, each component of this integration is ($\forall j\in[m]$)
\begin{align*}
     &\, \frac{1}{T}\int_{0}^T \!g_j(\bx+ \epsilon_a\bd(\bom t)) \, dt \\
= & \,\frac{1}{T}\int_{0}^T \!
g_j(\bx) + \epsilon_a  \nabla g_j(\bx)^\top \bd (\bom t) +\sO(\epsilon_a^2)\, dt\\
= & \,g_j(\bx) + \sO(\epsilon_a^2). 
\end{align*}
Combining these two parts, Lemma \ref{lemma:aveg2} is proved.

\subsection{Proof of Lemma \ref{lemma:omega}}\label{app:omega}


Take $\epsilon_a$ sufficiently small such that $\sO(\epsilon_a)< 1$. 
For any precision $\nu\in(0,1)$, there exists a time $T_1>t_0$ such that for any $t\geq T_1$,  $\beta(\Delta, t-t_0)\leq \frac{\nu}{4}$. Such  $T_1$ always exists because $\beta$ is a class-$\mathcal{KL}$ function, and thus  $||\bar{\bz}(t)||_{\hat{\sA}}\leq \frac{\nu}{2}$ for $t\geq T_1$ by \eqref{eq:y1con}. In addition,
by the exponential input-to-output stability of the fast dynamics in \eqref{eq:realave}, there exists a time $T_2> t_0$ such that for any $t\geq T_2$, every solution of \eqref{eq:realave} with $\bar{\bs}(t_0)\in [(\hat{\sA}\!+\!\Delta \mathbb{B})\cap \hat{\sZ}]\!\times\! \Delta \B$ satisfies 
\begin{align}
    ||\bar{\psi}(t)||\leq \frac{\nu}{2} + \sup_{\tau\geq t_0} ||\bm{\ell}(\bar{\bz}(\tau)) +\sO(\epsilon_a)||\leq \frac{\nu}{2}+ M_2.
\end{align}
Thus, for all $t\geq \max\{T_1,T_2\}$, the trajectory $\bar{\bs}(t)$ converges to a $\frac{\nu}{2}$-neighborhood of $\hat{\sA}\times M_2\B$.
Since the Omega-limit set from $\sF\times M_3\B$ is nonempty and $\Omega(\sF\times M_3\B)\subset (\hat{\sA}\times M_2\B) + \frac{\nu}{2} \B \subset \mathrm{int}(\sF\times M_3\B)  $.  By \cite[Corollary 7.7]{Goebel:12}, the set $\Omega(\sF\times M_3\B)$ is uniformly globally asymptotically stable for the average system \eqref{eq:realave} restricted to $\sF\times M_3\B$. 

%% file: Reference.tex
\bibliography{IEEEabrv, ref}